\documentclass[11pt]{amsart}

\usepackage{a4wide,amsmath,amssymb}
\usepackage{dsfont}
\usepackage{epsfig}
\usepackage{graphicx}
\usepackage{subfig}
\usepackage{verbatim}

\newtheorem{theorem}{Theorem}

\newtheorem{lemma}[theorem]{Lemma}

\newcommand{\ts}{\hspace{0.5pt}}
\newcommand{\C}{\mathbb{C}\ts}

\newcommand{\K}{\mathbb{K}\ts}

\newcommand{\Q}{\mathbb{Q}\ts}
\newcommand{\R}{\mathbb{R}\ts}

\newcommand{\one}{\mathds{1}}

\newcommand{\N}{\mathbb{N}}
\newcommand{\Z}{\mathbb{Z}}

\newcommand{\cE}{\mathcal{E}}

\def\mod#1{\,({\rm mod\ }#1)}
\def\proof{\noindent{\sc Proof.}\hskip 5pt}
\def\endproof{\hfill\vbox{\hrule
    \hbox{\vrule\kern4pt\vbox{\kern4pt
    \kern4pt}\kern4pt\vrule}\hrule}\bigskip}

\DeclareMathOperator{\GL}{GL}

\DeclareMathOperator{\trace}{tr}

\begin{document}

\title[Heights in piecewise affine maps]
{Growth of heights\\ in piecewise-affine planar maps}

\author{John A.~G.~Roberts}
\address{School of Mathematics and Statistics,
University of New South Wales,
Sydney, NSW 2052, Australia}
\email{jag.roberts@unsw.edu.au}
\urladdr{http://www.maths.unsw.edu.au/\~{}jagr}

\author{Franco Vivaldi}
\address{School of Mathematical Sciences, Queen Mary,
University of London,
London E1 4NS, UK}
\email{f.vivaldi@maths.qmul.ac.uk}
\urladdr{http://www.maths.qmul.ac.uk/\~{}fv}

\begin{abstract}
We consider the growth of heights of the points of the orbits of 
(piecewise) affine maps of the plane, with rational parameters.
We analyse the asymptotic growth rate of both global and local ($p$-adic) 
heights, for the primes $p$ that divide the parameters.
We show that almost all the points in a domain of linearity (such as an elliptic 
island in an area-preserving map) have the same exponential growth rate. 
We also show that the convergence of the $p$-adic height may be non-uniform, 
with arbitrarily large fluctuations occurring arbitrarily close to any point.
We explore numerically the behaviour of heights in the chaotic regions,
in both area-preserving and dissipative systems.
\end{abstract}
\date{\today}

\maketitle

\section{Introduction}\label{section:Introduction}

This paper is concerned with the growth rate of some indicators of 
arithmetical complexity ---the global and local (or $p$-adic) 
heights--- of the points of the orbits of affine and piecewise affine 
planar maps.
We present a combination of rigorous results and numerical experiments
connecting growth of heights to the dynamics on a divided phase space, 
where regular and irregular motions co-exist (see figure 1). 
This programme aims to develop a local analogue of the
so-called integrability criteria, which are detectors of global 
regularity of motions. These criteria have been the object of 
extended investigations; in particular, the notion of
diophantine integrability has been recently suggested, which 
is based on the slow growth of global heights ---see
\cite{Halburd} and references therein.
\begin{figure}[h]
\epsfig{file=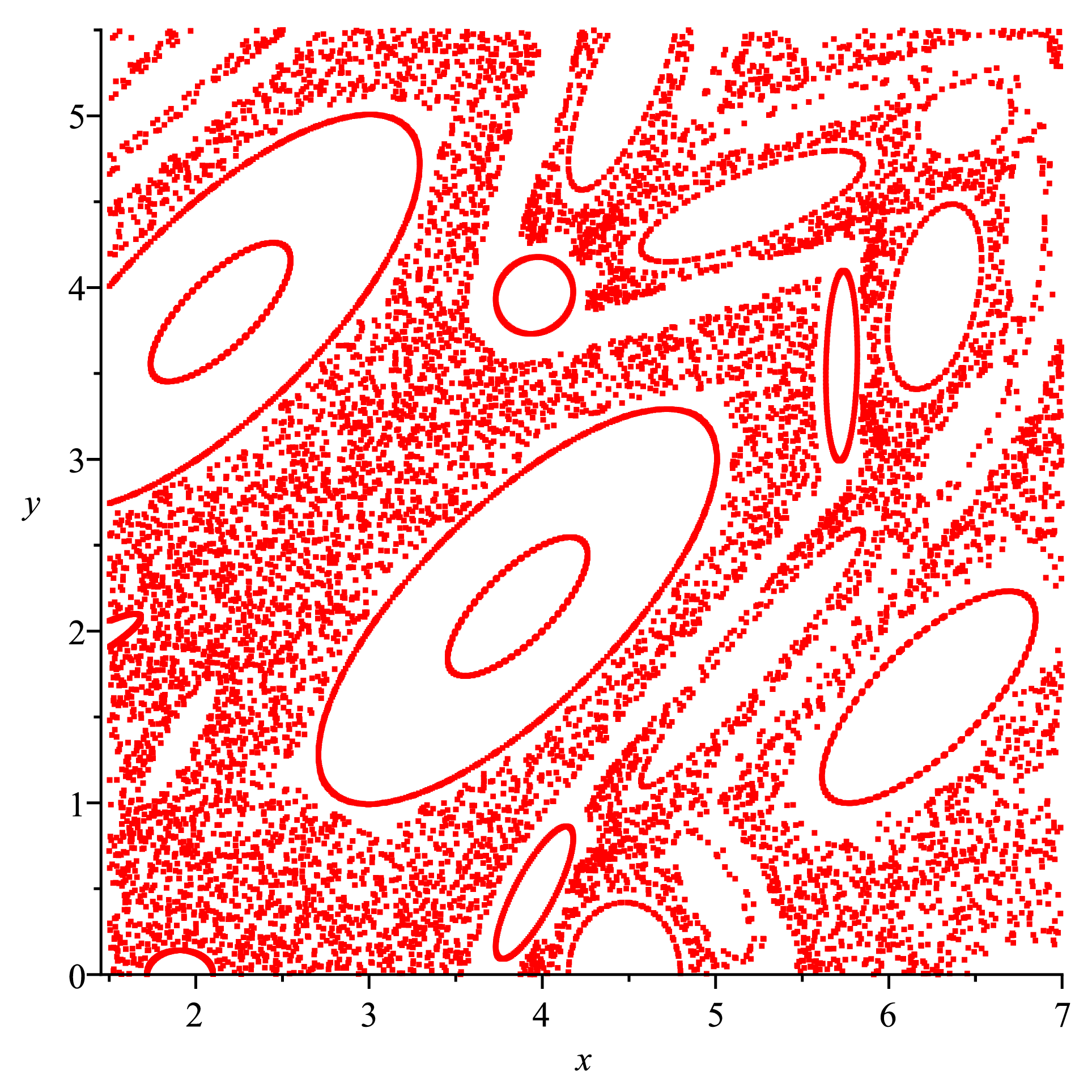,width=8cm,height=8cm}
\vspace*{-5pt}
\caption{\label{fig:PhasePortrait}\rm\small
Phase portrait of the area-preserving map $F$ defined in equation 
(\ref{eq:F}), with $f$ given in (\ref{eq:Myf}) and $d=1$, showing
a mixture of regular orbits on island chains and chaotic 
orbits.
}
\end{figure}

We are interested in monitoring the arithmetical complexity of the
points of an orbit of a piecewise affine map $F:\Q^2\to\Q^2$.
The simplest measure of the complexity of a rational number $x=m/n$ is its 
{\it height} \/ $H(x)$, defined as \cite[chapter 3]{Silverman}
\begin{equation}\label{eq:Height}
H(m/n)=\max(|m|,|n|)\hskip 40pt
\gcd(m,n)=1.
\end{equation}
The notions of size and height are extended to two dimensions as follows
\begin{equation}\label{eq:NormHeight}
\Vert z\Vert=\max(|x|,|y|)
\hskip 30pt
H(z)=\max(H(x),H(y))\qquad z=(x,y).
\end{equation}
The height will typically grow exponentially along orbits, so we define 
an allied quantity, the {\it logarithmic height}:
\begin{equation}\label{eq:LogHeight}
h(z)=\lim_{t\to\infty} \frac{1}{t}\,\log H(F^t(z))
\end{equation}
if the limit exists.
We have $h(z)=h(F(z))$, so the logarithmic height is a property of an orbit.
If $z$ is a (pre)-periodic point, then $H(F^t(z))$ is bounded, so that 
$h(z)=0$ (as long as the orbit of $z$ doesn't go through the origin).

Further indicators of complexity are defined by means of the $p$-adic absolute 
value $|\,\cdot\,|_p$, where $p$ is a prime number. 
(For background reference on $p$-adic numbers, see \cite{Gouvea}.)
Let the {\it order} $\nu_p(m)$ of an integer $m$ be the 
largest non-negative integer $k$ such that $p^k$ divides $m$, with $\nu(0)=\infty$. 
This definition is extended to the rational numbers $r=m/n$ by letting
$\nu_p(r)=\nu_p(m)-\nu_p(n)$ (the value of this expression doesn't depend 
on $m$ and $n$ being co-prime). Finally, we define
$$
|\,r\,|_p= p^{-\nu_p(r)}.
$$
The function $|\,\cdot\,|_p:\Q\to \Q$ has the properties of the
ordinary absolute value, with the triangular inequality replaced by the
stronger ultrametric inequality
\begin{equation}\label{eq:Ultrametric}
|x+y|_p\leqslant \mathrm{max}(|x|_p,|y|_p)
\qquad\mbox{or}\qquad
\nu_p(x+y)\geqslant \mathrm{min}(\nu_p(x),\nu_p(y))
\end{equation}
where equality holds if $|x|_p\not=|y|_p$ (or $\nu_p(x)\not=\nu_p(y)$).
We shall be using the estimate
\begin{equation}\label{eq:ValuationEstimate}
\nu_p(n) \leqslant \frac{\log n}{\log p} \qquad n\geqslant 1.
\end{equation}

The following identity connects the various absolute values over $\mathbb{Q}$:
\begin{equation}\label{eq:Adelic}
\forall x\in\Q\setminus \{0\},\qquad |x|\,\prod_p |x|_p =1
\end{equation}
where the product is taken over all primes.
Only finitely many terms of this product are different from 1; they correspond
to the prime divisors of the numerator and the denominator of $x$.

In two dimensions we use the quantities
\begin{equation}\label{eq:p-NormHeight}
\Vert z\Vert_p=\max(|x|_p,|y|_p)
\hskip 30pt 
\nu_p(z)=\min(\nu_p(x),\nu_p(y)).
\end{equation}
The norm $\Vert \cdot\Vert_p$ and valuation $\nu_p$ can be shown to satisfy the 
ultrametric inequalities analogous to (\ref{eq:Ultrametric}), respectively, with
equality holding if the two terms have distinct size.
Next we define the analogue of (\ref{eq:LogHeight}), namely the
\textit{$p$-adic} (or \textit{local}) \textit{height} $h_p(z)$ of the
initial point $z$ of an orbit:
\begin{equation}\label{eq:LocalHeight}
h_p(z)=\lim_{t\to\infty}\, -\frac{1}{t}\,\nu_p(F^t(z)).
\end{equation}
Comparing (\ref{eq:LocalHeight}) with (\ref{eq:LogHeight}), we note that the function 
$\nu_p$ is already logarithmic, and that there is no need of considering separately 
numerator and denominator, since the prime $p$ will appear only in one of them.

The functions $h$ and $h_p$ are variants of the so-called {\it canonical height} 
defined for rational functions of degree greater than one \cite[chapter 3]{Silverman}.
In this case, in place of (\ref{eq:Height}) one defines
$$
\hat H(m/n)=\max(|m|,|n|)\prod_p\max(|m|_p,|n|_p)
$$
and then one lets
$$
\hat h(x)=\lim_{t\to\infty}\frac{1}{\mbox{\rm deg}(F)^t}\log H(F^t(x))
$$
where $\mbox{\rm deg}(F)>1$ is the degree of $F$.
The height $\hat h$ behaves nicely under iteration: $\hat h(F(x))=\mbox{deg}(F)\hat h(x)$.
It measures the average rate of growth of the degree of $F$, collecting 
contributions from all absolute values.

In the case of (piecewise) affine mappings, the increase in complexity does not 
derive from degree growth, but rather from the growth of coefficients, hence
the definition of $h$ and $h_p$. 
Furthermore, we have kept the contributions from the various primes separate 
(as in the so-called {\it local canonical heights}) because they contain 
valuable information about the dynamics.

The height may be used to characterize generic properties of rational points.
To this end, we consider the set $\mathcal{B}_N$ of points in $\Q^2$ 
whose height is at most $N$:
\begin{equation}\label{eq:BoundedHeight}
\mathcal{B}_N=\{z\in\Q^2\,:\,H(z)\leqslant N\}.
\end{equation}
This set is finite. Indeed if $H(m/n)\leqslant N$, then 
$H(-m/n),H(\pm n/m)\leqslant N$, and we deduce that
$$
\#\mathcal{B}_N= \bigl(3+4\sum_{k=2}^N\phi(k)\bigr)^2\sim \frac{12^2}{\pi^4} N^4\qquad (N\to\infty)
$$
where $\phi$ is Euler's function \cite[section 5.5]{HardyWright} 
and where we have used the estimate $\sum_{k=1}^N\phi(k)\sim 3N^2/\pi^2$ 
(see \cite[theorem 330]{HardyWright} and also \cite[p 135]{Silverman}).
Half of the elements of $B_N$ lie within the square $\Vert z \Vert\leqslant 1$, 
where they approach a uniform distribution (because the Farey sequence has that 
property \cite{Mikolas,CodecPerelli}); the other half lie outside the square, and
they are obtained from the points inside the square by an inversion.
Thus the limiting distribution of points of bounded height approaches a smooth
limit on sufficiently regular bounded sets.

Let us now consider a set $A$ such that $A\subset X\subset \Q^2$, where
$X$ is some ambient set (possibly the whole of $\mathbb{Q}^2$). 
The density $\mu(A)$ of $A$ (in $X$) with respect to $\mathcal{B}_N$ is given by
\begin{equation}\label{eq:Density}
\mu(A)=\lim_{N\to\infty}\frac{\#(A\cap\mathcal{B}_N)}{\#(X\cap\mathcal{B}_N)}
\end{equation}
if the limit exists\footnote{For this it suffices to require that the closure of the
boundary of $A$ has zero measure (Jordan measurability)}. 
If $\mu(A)=1$, then we say that $A$ is `generic', or that the defining 
property of $A$ holds `almost everywhere' (in $X$). For example, the rational 
points on a smooth curve on the plane have zero density and hence are non-generic.

For the numerical experiments reported in section \ref{section:NumericalExperiments} 
we have chosen maps $F$ of the form
\begin{equation}\label{eq:F}
F:\R^2\to\R^2
\qquad
(x,y)\mapsto (f(x)-y,dx)
\end{equation}
where $f$ is a piecewise-affine real function and $d$ is a real number (the
Jacobian determinant of $F$). More precisely, we have a set $I$ of indices (possibly 
infinite), a partition $\{\Delta_i\}_{i\in I}$ of the real 
line into intervals, and a collection $\{f_i\}_{i\in I}$ of real affine 
functions 
$$
f_i:\R\to\R
\qquad
x\mapsto a_ix+b_i
\qquad
a_i, b_i\in\R
$$
such that
$$
f(x)=f_i(x)\qquad x\in \Delta_i.
$$
If $d=1$, then for any choice of $f$ the map $F$ is area-preserving
(see section \ref{section:PiecewiseAffineMaps}).
The literature devoted to maps of this type is substantial 
\cite{Devaney,BirdVivaldi,BeardonBullettRippon,AharonovDevaneyElias,%
LagariasRains,LagariasRains:b,LagariasRains:c}.

Let now $a_i,b_i,d\in\Q$. 
Then the set $\Q^2$ is invariant under $F$, and it makes sense to 
restrict the dynamics to rational points. (In fact one can restrict the 
space further ---see the appendix.)

\begin{figure}[t]
\hfil
\epsfig{file=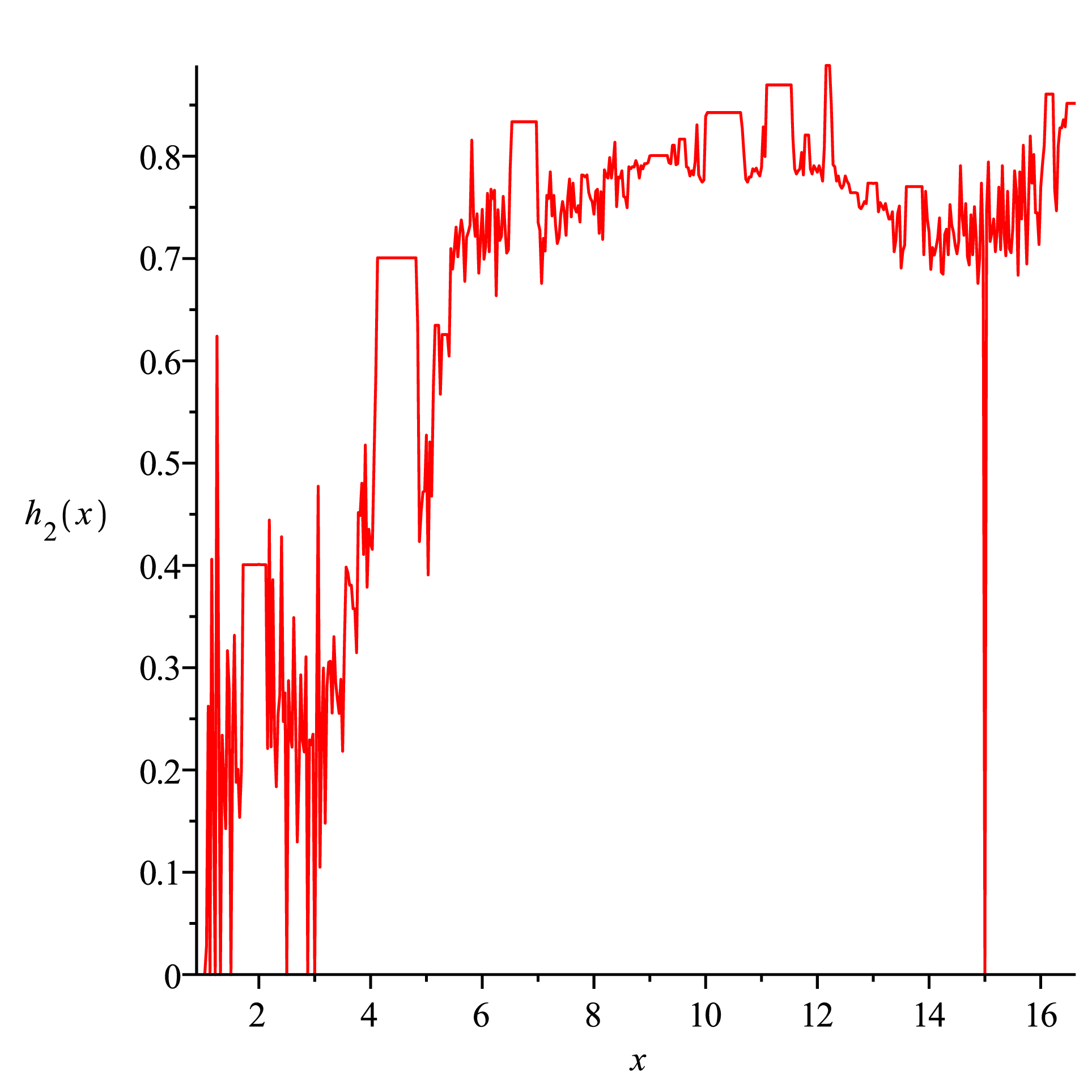,width=7cm,height=7cm}
\quad
\epsfig{file=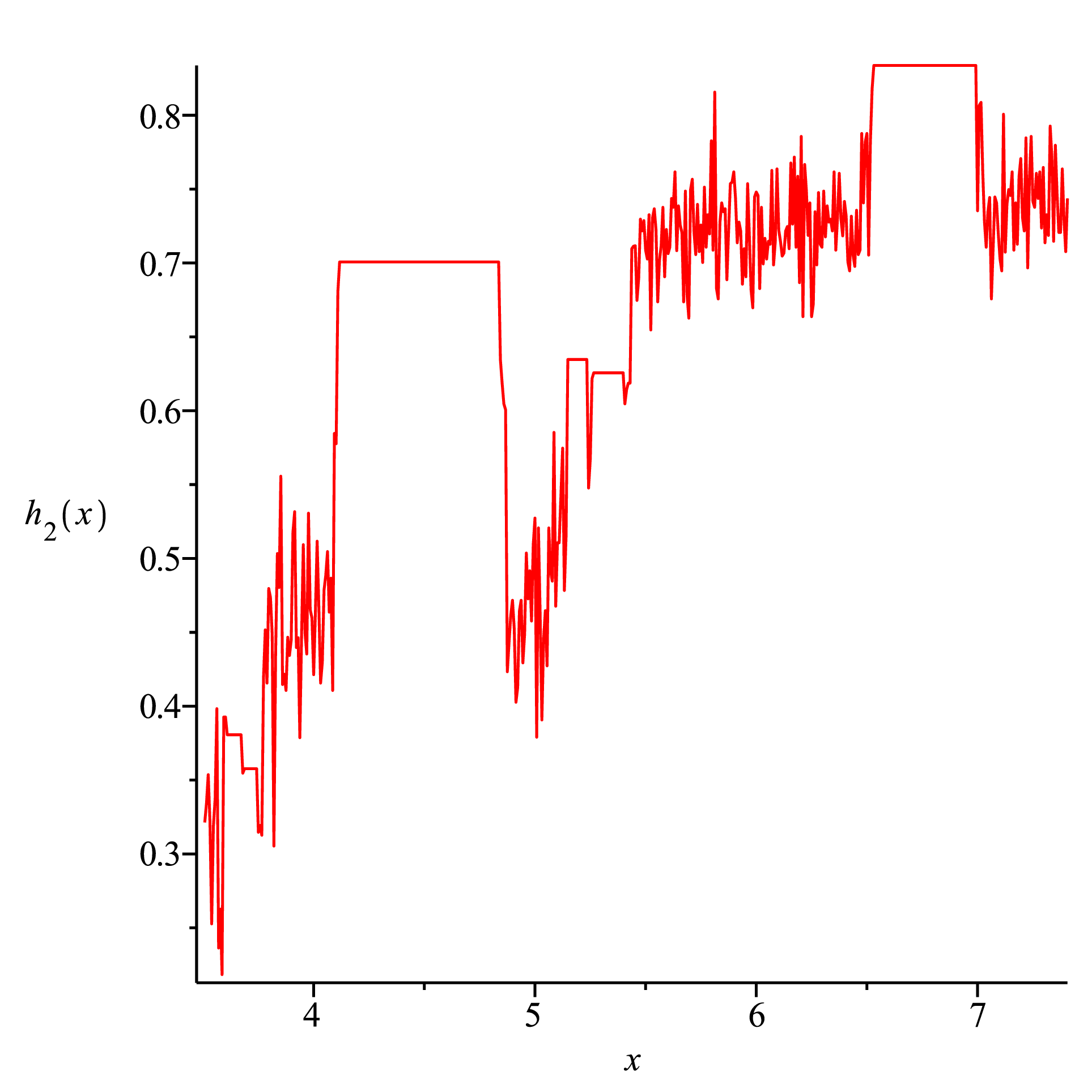,width=7cm,height=7cm}\hfil
\vspace*{-5pt}
\caption{\label{fig:Height}\rm\small
Behaviour of $h_2(x)$ for the map $F$ defined in equation 
(\ref{eq:Myf}), with initial conditions $z_0=(x,0)$. The plot
on the right shows a detail of that on the left.
}
\end{figure}

The $2$-adic height for some orbits of the map $F$ given by 
\begin{equation}\label{eq:Myf}
f(x)=\begin{cases}
 \frac{3}{2} x+\frac{3}{2} & x < -1\\
 0                         & -1 \leqslant x \leqslant  1\\
 \frac{3}{2} x-\frac{3}{2} & x > 1
\end{cases}
\end{equation}
with $d=1$ is shown in figure \ref{fig:Height}. 
The initial conditions are evenly spaced rational points on the positive $x$-axis.
The alternation of constancy and fluctuations is a distinctive feature of
height functions along smooth curves in phase space, which is connected
to the co-existence of regular and irregular motions. 
(To wit, compare figures \ref{fig:PhasePortrait} and \ref{fig:Height}.)

The plan of this paper is the following.
In section \ref{section:LocalHeights} we compute the local height
in affine maps, and show that, generically, all rational
points have the same height (theorem \ref{thm:LocalHeight}).
We identify the conditions under which convergence of the heights is 
non-uniform, but also show that the set of points having slow convergence 
have an exponentially large global height.
We then obtain explicit formulae for the valuation function $\nu_p$ along 
orbits in terms of Lucas polynomials; this gives us an alternative proof 
of theorem \ref{thm:LocalHeight}.
In section \ref{section:GlobalHeight} we determine the global height of an 
affine map, and show that, generically, all rational points have the same 
height (theorem \ref{thm:LogH}).
In section \ref{section:PiecewiseAffineMaps} we consider piecewise-affine maps
$F$ of $\Q^2$ (which include maps of the type (\ref{eq:F})), and their islands, 
which are bounded invariant domains where the motion is locally linear.
In the islands the results of the previous sections apply and all heights are
constants, which explains the plateaus in figure \ref{fig:Height} 
(theorem \ref{theorem:ConstantHeight}). 

In section \ref{section:NumericalExperiments} we explore numerically the convergence 
of height functions in the chaotic regions, and also consider briefly heights of 
quasi-periodic points.
In the appendix we construct a set $\mathbb{L}^2$, where $\mathbb{L}$ is a module
over a certain sub-ring of $\Q$ depending on the map's parameters,
which serves as a natural minimal phase space of a piecewise affine map.
This is the set relevant to our numerical experiments. 

\bigskip\noindent
{\sc Acknowledgements:} \/ 
JAGR and FV would like to thank, respectively, the School of 
Mathematical Sciences at Queen Mary, University of London, 
and the School of Mathematics and Statistics at the University 
of New South Wales, Sydney, for their hospitality.
This work was supported by the Australian Research Council.

\section{Local heights in affine maps}\label{section:LocalHeights}

We consider the behaviour of local heights (\ref{eq:LocalHeight}) of the
rational points for the affine map:
\begin{equation}\label{eq:Faffine}
F:\Q^2\to\Q^2
\qquad
z=(x,y) \mapsto \mathrm{M}\, z + s
\end{equation}
where $\mathrm{M}\in\GL(2,\Q)$ is a non-singular matrix with rational entries, 
and $s\in\Q^2$. (For notational ease, we do not use transpose symbols 
where it is clear by context, e.g., for $z$ and $s$ above.)

The map $F$ has a single rational fixed point 
\begin{equation}
z^*=(x^*,y^*)=-(\mathrm{M}-\one)^{-1}\,s,\nonumber
\end{equation}
and if $z_0=z^*+z_0^\prime$, then
\begin{equation} \label{eq:Iterate}
 z_t=F^t(z_0)=\mathrm{M}^t\,z_0^\prime+z^*.
\end{equation}
We define
\begin{equation} \label{eq:TD}
T=\trace (\mathrm{M}), \hskip 30pt D=\det (\mathrm{M}),
\end{equation}
and we let $q(x)=x^2-Tx+D$ be the characteristic polynomial of $\mathrm{M}$,
with roots $\alpha$ and $\beta$. 

\begin{figure}[t]
\hfil
\epsfig{file=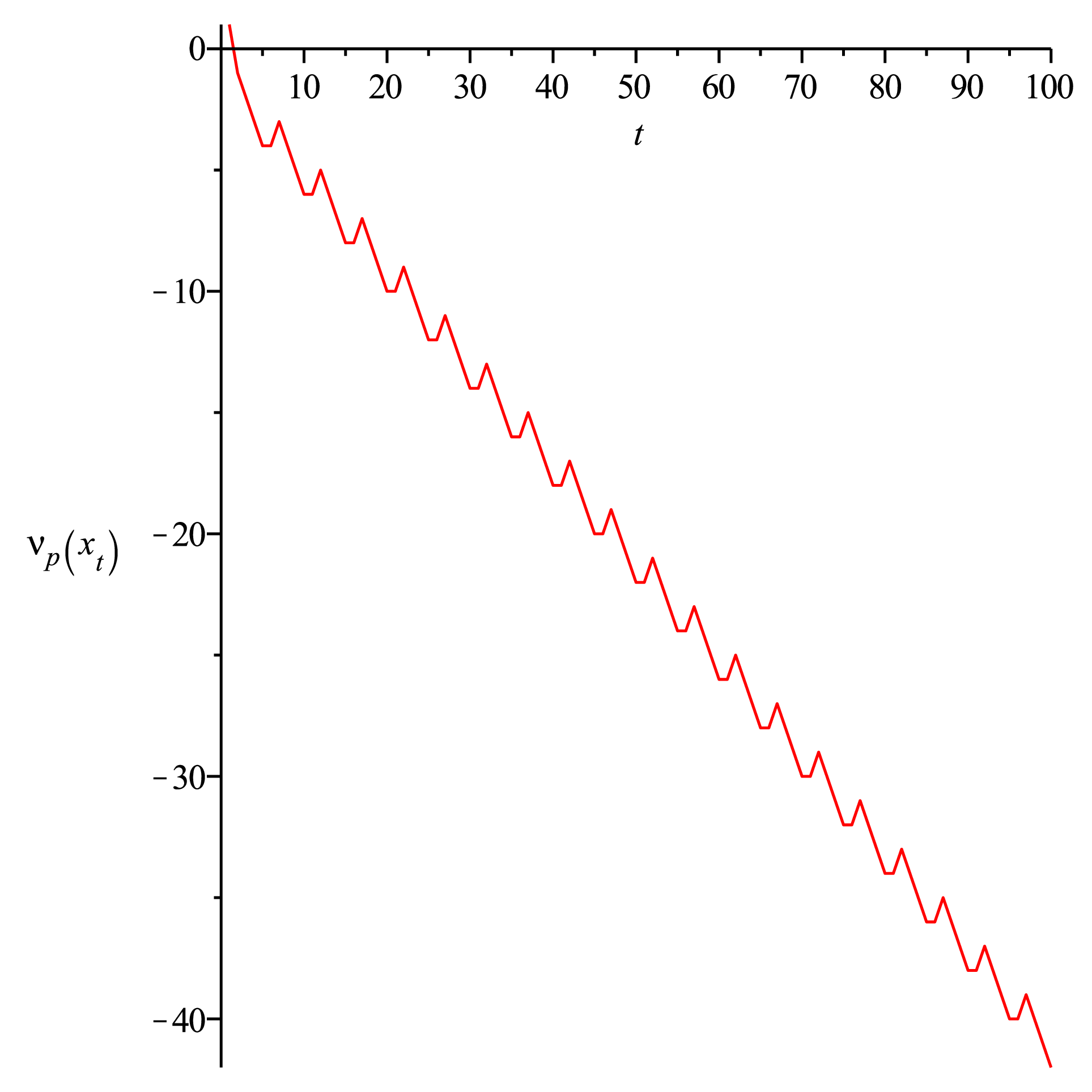,width=7cm,height=6cm}
\quad
\epsfig{file=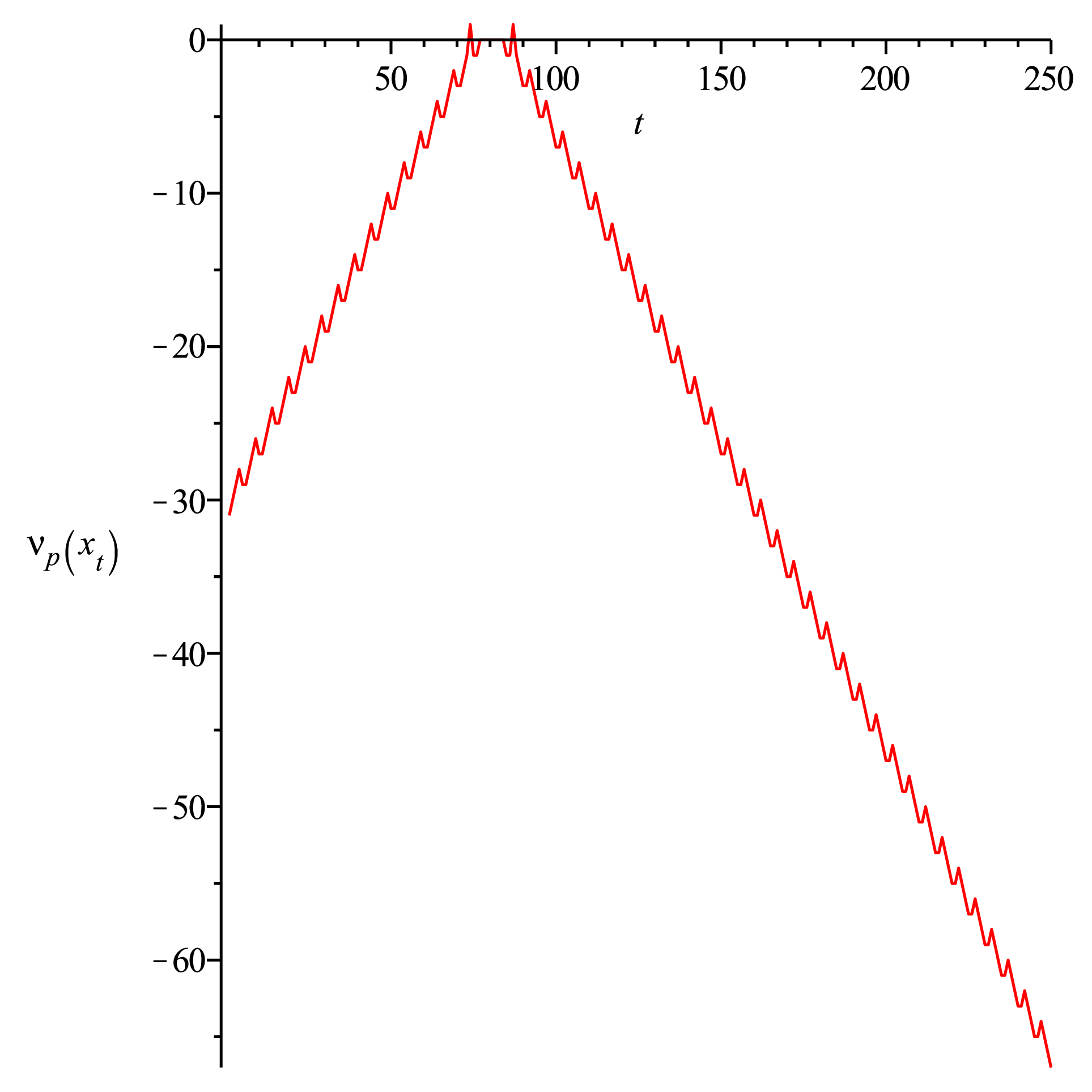,width=7cm,height=6cm}\hfil
\vspace*{-5pt}
\caption{\label{fig:nu}\rm\small
Time-dependence of $\nu_2(x_t)$ for two rational orbits 
of the map (\ref{eq:Myf}), with very close initial conditions inside 
the same island with elliptic periodic point $z^*=(21/11,0)$. 
Left: typical behaviour, for $z_0=(2,0)$. Right: anomalous
behaviour, for $z_0=(2,0)+z'$ with $\Vert z'\Vert<10^{-8}$.
In this case the point $z_0$ lies in the vicinity of the stable manifold 
of $z^*$ in $\Q_2^2$.
}
\end{figure}

The computation of $p$-adic heights is an eigenvalue problem analogous to the 
computation of the Lyapunov exponent. 
Further insight is obtained by studying the detailed behaviour of the sequence
$(\nu_p(z_t))$ (see figure \ref{fig:nu}), which will be considered in 
section \ref{section:ExplicitFormulae}.

\begin{theorem}\label{thm:LocalHeight}
Let $F$ be the affine map (\ref{eq:Faffine}) with $T,D,s$ as above.
If $s=(0,0)$, then for almost all $z\in\Q^2$ we have:
\begin{enumerate}
\item [$i)$] if $\nu_p(D)>2\nu_p(T)$ then $h_p(z)=-\nu_p(T)$;
\item [$ii)$] if $\nu_p(D)\leqslant 2\nu_p(T)$ then $h_p(z)=-\nu_p(D)/2$.
\end{enumerate}
If $s\not=(0,0)$ then the above expressions for $h_p$ must be replaced by 
$\mathrm{max}(-\nu_p(T),0)$ and $\mathrm{max}(-\nu_p(D)/2,0)$, respectively.
\end{theorem}

\smallskip\noindent {\sc Proof.}
Let $\mathbb{Q}_p$ be the completion of $\mathbb{Q}$ with respect to the absolute 
value $|\cdot|_p$. The eigenvalues $\alpha,\beta$ of $\mathrm{M}$ lie in a field
$K$ which is either $\mathbb{Q}_p$ or a quadratic extension of $\mathbb{Q}_p$. 
In $K$ there is a prime element $\pi$ (either $p$ or $\sqrt{p}$) and a valuation 
$\nu_\pi$, which is either $\nu_p$ or is an extension of $\nu_p$ which agrees 
with $\nu_p$ on $\mathbb{Q}$. 
Let $\alpha$ be a largest eigenvalue, that is,
$\nu_\pi(\alpha)\leqslant\nu_\pi(\beta)$. 
Let 
\begin{equation}\label{eq:uv}
u=\nu_p(D),\qquad v=\nu_p(T)
\end{equation}
and let $\Pi$ be the Newton polygon of $q(x)$, namely the convex hull of 
the points $(0,\infty)$, $(0,u)$, $(1,v)$, $(2,0)$, $(2,\infty)$. 
If $u>2v$ then $\Pi$ has two finite sides with distinct slopes 
$v-u$ and $-v$, of which the latter is the largest. Hence by
\cite[Theorem 6.4.7]{Gouvea} we have
$\nu_\pi(\beta)=u-v$ and $\nu_\pi(\alpha)=v$.
Likewise, if $u\leqslant 2v$ then $\Pi$ has one finite side of slope $-u/2$.
Hence $\nu_\pi(\beta)=\nu_\pi(\alpha)=u/2$. 

First we consider the parameter $s=(0,0)$. 
We begin with the case $|\alpha|_\pi>|\beta|_\pi$, which is case $i)$. 
We have $K=\mathbb{Q}_p$ (see section \ref{section:Eigenspaces}), 
and we write 
\begin{equation}\label{eq:LinearSolution}
z_t=\alpha^t c_1\mathbf{w}_1+\beta^t c_2\mathbf{w}_2
\end{equation}
where the $\mathbf{w}_i$ are linearly independent eigenvectors of $\mathrm{M}$ 
in $\mathbb{Q}_p^2$ and the coefficients $c_i$ are in $\mathbb{Q}_p$.
For generic initial conditions $c_1\not=0$ (i.e., $z_0$ does not lie in the eigenspace generated 
by $\mathbf{w}_2$), the $p$-adic height in a linear system is determined by the 
eigenvalue with largest $p$-adic absolute value, which is $\alpha$. 
Specifically, for all large enough $t$, the two terms in (\ref{eq:LinearSolution}) 
have distinct size, and hence from (\ref{eq:p-NormHeight}) and following comments we see that 
\begin{equation}\label{eq:NormOfz_t}
\Vert z_t\Vert_p=|\alpha|_p^t\Vert c_1\mathbf{w}_1\Vert_p,
\end{equation}
from which $\nu_p(z_t)\sim t\nu_p(\alpha)$ and the result follows.

Let us now deal with case $ii)$.
If $|\alpha|_\pi=|\beta|_\pi$, but $\alpha\not=\beta$, we rewrite (\ref{eq:LinearSolution}) 
as
$$
z_t=\alpha^t u_t
\hskip 40pt u_t=c_1\mathbf{w}_1+(\beta/\alpha)^t c_2\mathbf{w}_2,
$$
noting that $\alpha$ is non-zero. Then $\mu=\beta/\alpha$ is a $p$-adic unit, 
and hence there exits a smallest positive integer $n$ such that
$\mu^n=\overline\mu=1+\gamma$ with $|\gamma|_\pi<1$. If $\gamma=0$, that is,
$\mu$ is a root of unity, then $u_t$ is periodic, and hence 
$h_p(z)=-\nu_p(\alpha)$, as desired. 

If $\gamma\not=0$, then the sequence $(\overline{\mu}^t)$ is dense in a disc
(see \cite{ArrowsmithVivaldi} and \cite[chapter 5]{Hasse}), and hence $(\mu_t)$ is 
dense in the union of $n$ discs.
Thus each component of $z_t=(x_t,y_t)$ is also dense in a finite union of discs. 
If none of these discs contains the origin, then $\Vert u_t\Vert_p$ assumes 
finitely many values, and the result follows. Otherwise $\Vert u_t\Vert_p$ is 
bounded above but not bounded away from zero, and the rate at which $\Vert u_t\Vert_\pi$ 
approaches zero is the same as the rate at which $\overline \mu^t$ approaches 1.
From the binomial theorem we obtain $|\overline \mu^t-1|_\pi=p^{\nu_\pi(t)}$
and hence the quantity $\max_{t<T}\{\nu_\pi(z_t')\}$ grows logarithmically,
from (\ref{eq:ValuationEstimate}).
It follows that $\nu_p(z_t)\sim t\nu_\pi(\alpha)$, as desired.

Finally, if the Jordan form of $M$ is not diagonal, then the sequence $(z_t)$
contains a term affine in $t$. The contribution of this term is logarithmic,
again due to (\ref{eq:ValuationEstimate}). Hence, in all cases, $h_p=-\nu_\pi(\alpha)$. 

If $s\not=(0,0)$ then from (\ref{eq:p-NormHeight}) and (\ref{eq:Iterate})
we find $\nu_p(z_t)\geqslant \mathrm{min}(\nu_p(M^tz_0'),\nu_p(z^*))$.
In case $i)$, if $v<0$, then, for all sufficiently large $t$ the first term
is the largest, that is, $\nu_p(z_t)=\nu_p(M^tz_0')$, and the previous analysis
applies. Likewise, if $v>0$, then eventually the second term becomes the 
largest, and since this term is constant, we get $h_p=0$. If $v=0$, then
the inequality remains such, but the first term grows at most 
logarithmically, and so $h_p=0$. Case $ii)$ is treated similarly.
\hfill $\Box$

\subsection{$p$-adic eigenspaces}\label{section:Eigenspaces}

We look more closely at the $p$-adic dynamics of a linear map with eigenvalues
$\alpha,\beta$ of distinct magnitude, which is case $i$) of theorem \ref{thm:LocalHeight}.
Using the notation (\ref{eq:uv}), we see that if $u>2v$, then, necessarily, 
$v\not=+\infty$ ($T\not=0$).
Letting
$$
T=T'p^{\nu(T)}
\hskip 30pt
D=D'p^{\nu(D)}
$$
we have $T'\not=0$. 
Let now $\theta=p^{-v}\lambda$, where $\theta$ is a root of the polynomial
\begin{equation}\label{eq:CharPoly}
s(x)=x^2-T'x+D'p^{u-2v}
\qquad \mathrm{with}\qquad
\frac{ds(x)}{dx}=2x-T'.
\end{equation}
We have the factorisation:
\begin{equation}
s(x)\equiv x(x-T')\mod{p}.\nonumber
\end{equation}
Then $s(x)$ has two distinct roots modulo $p$, congruent to $0$
and $T'$, respectively, and at these roots $s(x)$ is equal to
$\pm T'\not\equiv 0$ from (\ref{eq:CharPoly}).
From Hensel's lemma \cite[section 3.4]{Gouvea}, we have that $s(x)$ has two 
distinct roots in $\mathbb{Z}_p$, which we denote by $\alpha',\beta'$, of which
the largest, $\alpha'$, is a unit.
Hence $\nu_p(\alpha)=\nu_p(T)$, in agreement with theorem \ref{thm:LocalHeight}

Now, the polynomial $s(x)$ is irreducible over $\Q$ if and only if $q(x)$ is
irreducible, since their roots differ by a rational factor.
If $q(x)$ is reducible, then these eigenspaces have infinitely many rational 
points; if $q(x)$ is irreducible, then these eigenspaces have no rational points, 
apart from the origin. 

In the first case there will be a non-generic (zero-density) set of rational 
points with height $\nu_p(T)-\nu_p(D)$, lying on the eigenspace corresponding to 
the smallest eigenvalue.
Thus a sufficient condition for all non-zero rational points to have the 
same $p$-adic height is $1=v\leqslant u$, for in this case $q(x)$ is 
irreducible by Eisenstein's criterion \cite[Proposition 5.3.11]{Gouvea}.

In the second case all points have the same height $-\nu_p(D)$, apart from the origin.
Rational approximants for the roots of $q(x)$ may be constructed by iterating Newton's 
map for $q(x)$ sufficiently many times, with an appropriate initial condition 
\cite[section 3.4]{Gouvea}. 
The components of an eigenvector of $\mathrm{M}$ may be chosen to be linear 
expression in such eigenvalues, with rational coefficients.

We are interested in motion in the $p$-adic vicinity of the eigenspace 
$W_p^\beta$ corresponding to the smaller eigenvalue. We begin with a general lemma.

\begin{lemma}\label{lemma:SimultaneousApproximation}
Let $p$ be a prime number. For any $z\in\Q^2$, any $\zeta\in\Q_p^2$,
and any $\epsilon>0$, there is $z'\in\Q^2$ such that
$$
\Vert z'-z\Vert +\Vert z'-\zeta\Vert_p<\epsilon
$$
with the norms (\ref{eq:NormHeight}) and (\ref{eq:p-NormHeight}), respectively.
\end{lemma}

\medskip

\proof
The rational sequence
\begin{equation}\label{eq:rk}
r_k=\frac{1}{1+p^k}\qquad k=1,2,\ldots
\end{equation}
has the property that, as $k\to\infty$, $r_k\to 0$ in $\Q$, while $r_k\to 1$ in $\Q_p$.
Let now $z=(x,y)\in\Q^2$ and $\epsilon>0$ be given. 
For any $(a,b)\in\Q^2$, the sequence 
\begin{equation}\label{eq:zk}
z^{(k)}=z+r_k(a,b)\qquad k=1,2,\ldots
\end{equation}
converges to $z$ in $\Vert\,\cdot\,\Vert$. We choose $K_1$ such that,
for all $k>K_1$, we have $\Vert z^{(k)}-z \Vert<\epsilon/2$.

Let $\zeta=(\zeta_1,\zeta_2)$. 
We will show that $a,b$ in (\ref{eq:zk}) may be chosen so that
$\Vert z^{(k)}-\zeta \Vert_p\to 0$.
We find
$$
z^{(k)}-\zeta=(x+ar_k-\zeta_1,\,y+br_k-\zeta_2).
$$
Since $\Q$ is dense in $\Q_p$, we can find $s=(s_1,s_2)\in\Q^2$ such that 
$\Vert \zeta-s\Vert_p<\epsilon/2$.
Let $a=s_1-x$. Then there is $K_2$ such that for all $k>K_2$
we have $|x+ar_k-s_{1}|_p<\epsilon/2$.
Similarly, let $b=s_2-y$. Then there is $K_3$ such that for all $k>K_3$
we have $|y+br_k-s_{3}|_p<\epsilon/2$.

Let now $K=\max(K_1,K_2,K_3)$. For all $k>K$, the ultrametric 
inequality (\ref{eq:Ultrametric}) gives
\begin{eqnarray*}
|x+ar_k-\zeta_1|_p 
 &=& |x+ar_k-s+s-\zeta_1|_p\\
 &\leqslant& \mathrm{max}(|x+ar_k-s|_p,|s-\zeta_1|_p)\\
 &\leqslant& \mathrm{max}(\epsilon/2,\epsilon/2)=\epsilon/2.
\end{eqnarray*}
Similarly, $|y+br_k-\zeta_2|_p \leqslant \epsilon/2$.
In the same $k$-range, we obtain
$$
\Vert z^{(k)}-\zeta\Vert_p=
\mathrm{max}(|x+ar_k-\zeta_1|_p,|y+br_k-\zeta_2|_p)
<\mathrm{max}(\frac{\epsilon}{2},\frac{\epsilon}{2})=\frac{\epsilon}{2}.
$$

We have shown that for all $k>K$, the point $z'=z^{(k)}$ lies
within an $\epsilon/2$-neighbourhood of $z$ in the ordinary 
norm, and within an $\epsilon/2$-neighbourhood of 
$\zeta$ in the $p$-adic norm. The lemma follows.
\endproof

Now choose $\zeta\in W_p^\beta\subset\Q_p^2$. The lemma states that
arbitrarily close to any rational point we can find another rational 
point as close as we please to an eigenvector $\zeta$ of $\mathrm{M}$.
Thus, irrespective of the rationality of the eigenvalues, there always will be 
a dense set of initial conditions that are to close to the eigenspace $W_p^\beta$ 
to cause the second term in (\ref{eq:LinearSolution}) to dominate for small values of $t$.
For these orbits the convergence of $h_p$ will be slow.
The sequence $(\nu_p(z_t))$ will feature two distinct affine regimes, with slopes 
$\nu_p(T)-\nu_p(D)$ and $-\nu_p(T)$, respectively. If the slopes have
different sign and $z^*\not=(0,0)$, then these regimes may be separated by a third 
regime, determined by a constant lower bound ---see figure \ref{fig:nu}.

We want to justify the statement that the height of a `typical' rational point 
converges rapidly to its asymptotic value $-\nu_p(T)$, in apparent defiance 
of the pathologies exposed by lemma \ref{lemma:SimultaneousApproximation} above.
We will show that points for which the non-archimedean height has anomalous 
time-dependence must also have a large archimedean height.
For brevity, we consider only the linear case.

Let $\cE=\mathbb{Q}^2\setminus W_p^{\beta}$.
Then, in the regime in which equation (\ref{eq:NormOfz_t}) holds, we have that 
$\Vert z_{t+1}\Vert_p=|\alpha|_p\Vert z_{t}\Vert_p$.
Now we define the \textit{lag time} $\tau(z)$ to be the time at which this 
asymptotic regime sets in, namely,
\begin{equation}\label{eq:tau}
\tau:\cE\to\N
\qquad
\tau(z)=\min\{t\in\N\,:\,
\forall s\geqslant t,\,\, 
\Vert z_{s+1} \Vert_p =|\alpha|_p\Vert z_s\Vert_p\}.
\end{equation}
Because the eigenspace $W_p^\beta$ of $\beta$ has been excluded, the function $\tau$ is 
well-defined. The larger the value of $\tau(z)$, the slower the convergence 
of the $p$-adic height $h_p(z)$.

From equations (\ref{eq:LinearSolution}) and (\ref{eq:tau}) and the ultrametric
inequality, we find that
$$
\left|\frac{\alpha}{\beta}\right|^{\tau(z)}=\left|\frac{c_2(z)}{c_1(z)}\right|_p
   \frac{\Vert \mathbf{w}_2\Vert_p}{\Vert \mathbf{w}_1\Vert_p}\,\gamma(z)
$$
where the quantity $\gamma\in (|\beta/\alpha|_p, 1]$ ensures that $\tau$ is an
integer. Hence, as $\tau\to\infty$ we must have $|c_2/c_1|_p\to\infty$.
Now, for any non-zero rational number $r$ and 
any prime $p$, we have the estimate $H(r)\geqslant p^{|\nu_p(r)|}$.
Hence for large enough $\tau$ there is a constant $\kappa$ independent of $z$ such that
$$
\kappa\left|\frac{\alpha}{\beta}\right|^{\tau(z)} \leqslant 
\left|\frac{c_2(z)}{c_1(z)}\right|_p =p^{\nu_p(c_1(z)/c_2(z))}
\leqslant H(c_1'(z)/c_2'(z)),
$$
where $c_1'$ and $c_2'$ are any rational approximants of $c_1$ and $c_2$ such that
$\nu_p(c_1'/c_2')=\nu_p(c_1/c_2)$. 
Thus the height of the ratio of the coefficients of $z_t$ in the representation
(\ref{eq:LinearSolution}) grows at least exponentially in the lag time $\tau$.

\subsection{Explicit formulae}\label{section:ExplicitFormulae}
 
In this section we derive explicit formulae for $z_t$ and $\nu_p(z_t)$, which will 
give us an alternative, more direct proof of theorem \ref{thm:LocalHeight}, 
with the exclusion of some special cases.

From \eqref{eq:Iterate}, we need the powers of the rational matrix $\mathrm{M}$. 
Using the Cayley-Hamilton theorem, one proves by induction
(e.g., \cite[Lemma 1]{BaakeRobertsWeiss}) that 
for $t\in \Z$, the following relation holds
\begin{equation} \label{eq:CayHam}
    \mathrm{M}^t \, = \, U_t\, \mathrm{M} - D\, U_{t-1}\, \one,
\end{equation}
where the sequence of rational numbers $U_t=U_t(T,D)$ obeys the recursion
\begin{equation} \label{eq:RecU}
U_0=0,\quad U_1=1,\qquad
U_{t+1}({T},{D})\,=\,{T}\ts U_{t}({T},{D}) - {D}\ts U_{t-1}({T},{D}), 
\quad t\geqslant 1.
\end{equation}
If $T$ and $D$ are integers, then $U_t$ is an integer sequence, known as
the {\em Lucas sequence of the first kind}. 
In a slight abuse of notation, we will use the same symbol $U_t$ for our 
case of a rational sequence generated by \eqref{eq:RecU} because
many of the properties of Lucas sequences are independent of whether 
$T$ and $D$ are integers.
It follows by iteration of \eqref{eq:RecU} that  $U_t({T},{D})$ 
is a polynomial in ${T}$ and ${D}$ with integer coefficients.
Its general form \cite{Hoggatt} is
\begin{equation}\label{eq:U}
U_t({T},{D}) =  \sum_{k=0}^{\lfloor (t-1)/2 \rfloor}\,c_k^{(t)}\;T^{t-2k-1}\,(-D)^k
\end{equation}
where
\begin{equation}
c_k^{(t)}=\binom{t-k-1}{k}.\nonumber
\end{equation}
We note that
\begin{equation}\label{eq:BinomialEstimate}
0\leqslant \nu_p\left( \binom{n}{m}\right) \leqslant 
  \left\lfloor\frac{\log n}{\log p}\right\rfloor-\nu_p(m).
\end{equation}
From (\ref{eq:U}), we have that, for all $t\geqslant 1$:
\begin{enumerate}
\item [--] The polynomial $U_t(T,D^2)$ is homogeneous of degree $t-1$.
\item  [--] The leading term of $U_t$ is ${T}^{t-1}$ (i.e., $U_t$ is monic) while
the term of lowest total degree is ${(-D)}^{ \frac{t-1}{2} }$ if $t$ is odd and 
$\frac{t}{2}\,{T}\,{(-D)}^{\frac{t}{2}-1}$ if $t$ is even. 
\end{enumerate}
 
From \eqref{eq:Iterate} with \eqref{eq:CayHam}, we see that
\begin{equation} \label{eq:Formzt}
z_t = \begin{pmatrix} {x}_t \\  {y}_t \end{pmatrix} =
U_t({T},{D})\, 
\begin{pmatrix} {x}_1^\prime \\  {y}_1^\prime \end{pmatrix}
- {D}\, U_{t-1}({T},{D})\, 
\begin{pmatrix} {x}_0^\prime \\  {y}_0^\prime \end{pmatrix}
+ 
\begin{pmatrix} {x}^* \\  {y}^* \end{pmatrix}
\end{equation}
where 
$$
\begin{pmatrix} {x}_1^\prime \\  {y}_1^\prime \end{pmatrix}={\mathrm{M}}\, 
   \begin{pmatrix} {x}_0^\prime \\  {y}_0^\prime \end{pmatrix}. 
$$ 

Let us now consider the first component of $z_t$ in \eqref{eq:Formzt}.
Using \eqref{eq:U} we rewrite it as follows:
\begin{equation} \label{eq:xt}
x_t  = \mathcal{T}_t^{(1)}+\mathcal{T}_t^{(0)}+{x}^* 
\end{equation}
where
\begin{equation} \label{eq:T1and0}
\mathcal{T}_t^{(1)}={x}_1^\prime \, \sum_{i_1=0}^{\lfloor (t-1)/2 \rfloor}\,c_{i_1}^{(t)}\,T^{t-2i_1-1}\,(-D)^{i_1},
\qquad
\mathcal{T}_t^{(0)}={x}_0^\prime \sum_{i_0=1}^{\lfloor t/2 \rfloor}\,c_{i_{0}-1}^{(t-1)}\,T^{t-2i_0}\,(-D)^{i_0-1}.
\end{equation}
The greatest value of the summation indices is given by:
\begin{eqnarray*}
t \mbox{ odd}: \quad i_1^{max}:=& \lfloor (t-1)/2 \rfloor=(t-1)/2 \quad &i_0^{max}:=\lfloor t/2 \rfloor=(t-1)/2 \\
t \mbox{ even}:\quad  i_1^{max}:=& \lfloor (t-1)/2 \rfloor=t/2-1 \quad &i_0^{max}:=\lfloor t/2 \rfloor=t/2. 
\end{eqnarray*}
From (\ref{eq:xt}) and the ultrametric inequality (\ref{eq:Ultrametric}) it follows that
\begin{equation}\label{eq:Orderxt}
\nu_p(x_t)\geqslant 
   \min(\nu_p(\mathcal{T}_t^{(1)}),\nu_p(\mathcal{T}_t^{(0)}),\nu_p({x}^*) ).
\end{equation}
For the order of the first term, using \eqref{eq:T1and0} gives
\begin{equation}
\nu_p(\mathcal{T}_t^{(1)})
\geqslant \nu_p(x_1^\prime)+\min_{i_1}(\nu_p(c_{i_1}^{(t)})+{i_1}\,(\nu_p(D)-2\nu_p(T))+(t-1)\nu_p(T)).\nonumber
\end{equation}
We have three cases:
\begin{enumerate}
\item [$i)$] $\nu_p(D)>2\nu_p(T)$.
Using (\ref{eq:BinomialEstimate}), we see that the unique minimum is 
achieved at $i_1=0$ with $c_{0}^{(t)}=1$, giving
$$
\nu_p(\mathcal{T}_t^{(1)})=\nu_p(x_1^\prime)+(t-1)\nu_p(T).
$$
\item [$ii)$] $\nu_p(D)<2\nu_p(T)$.
The unique minimum is achieved at $i=i_1^{max}$, where $c_{i_1^{max}}^{(t)}$is equal to $1$ 
when $t$ is odd and to $t/2$ when $t$ is even. Thus
$$
\nu_p(\mathcal{T}_t^{(1)})=\nu_p(x_1^\prime)+
\begin{cases}
\displaystyle \frac{t-1}{2}\, \nu_p(D) & t\, \mbox{odd}\\
\displaystyle \frac{t-2}{2}\,  \nu_p(D) +\nu_p\bigl(\frac{t}{2}\bigr) + \nu_p(T) & t\, \mbox{even.}
         \end{cases}
$$
\item [$iii)$] $\nu_p(D)= 2\nu_p(T)$.
A minimum is achieved at $i_1=0$, with $c_{0}^{(t)}=1$, giving
$$
\nu_p(\mathcal{T}_t^{(1)})\geqslant \nu_p(x_1^\prime)+(t-1)\nu_p(T).
$$
\end{enumerate}

A very similar analysis for the order $\nu_p(\mathcal{T}_t^{(0)})$ in \eqref{eq:Orderxt} gives
\begin{enumerate}
\item [$i)$] $\nu_p(D)>2\nu_p(T)$.
$$
\nu_p(\mathcal{T}_t^{(0)})=\nu_p(x_0^\prime)+(t-2)\, \nu_p(T)+ \nu_p(D).
$$
\item [$ii)$] $\nu_p(D)<2\nu_p(T)$.
$$
\nu_p(\mathcal{T}_t^{(0)})=\nu_p(x_0^\prime)+
\begin{cases}
\displaystyle  \frac{t-1}{2}\, \nu_p(D) + \nu_p\bigl(\frac{t-1}{2}\bigr) +\nu_p(T) & t\, \mbox{odd}\\
\displaystyle  \frac{t}{2}\,  \nu_p(D) & t\, \mbox{even.}
         \end{cases}
$$
\item [$iii)$] $\nu_p(D)= 2\nu_p(T)$.
$$
\nu_p(\mathcal{T}_t^{(0)})\geqslant \nu_p(x_0^\prime)+t \nu_p(T).
$$
\end{enumerate}

The analysis for the second component $y_t$ in (\ref{eq:Formzt}) is identical.
 
From the above and (\ref{eq:Orderxt}), we have:

\bigskip\noindent
$i)$\/  $\nu_p(D)>2\nu_p(T)$:
\begin{equation}\label{eq:InequalityI}
\nu_p(x_t)\geqslant \min\{\nu_p(x_1^\prime)+(t-1)\nu_p(T),
      \nu_p(x_0^\prime)+(t-2)\,\nu_p(T)+ \nu_p(D),
      \nu_p(x^*)\}.
\end{equation}
If $\nu_p(T)\geqslant 0$, then the linear terms are increasing, and we have two possibilities.
If $x^*\not=0$, then eventually we have $\nu_p(x_t)=\nu_p(x^*)$. 
If, $x^*=0$, then eventually, under the non-degeneracy condition
\begin{equation}\label{eq:NonDegeneracy}
\nu_p(x_1')+\nu_p(T)\not=\nu_p(x_0')+\nu_p(D)
\end{equation}
a unique minimum emerges in (\ref{eq:InequalityI}), and $\nu_p(x_t)$ 
becomes affine. In the degenerate case, the inequality
(\ref{eq:InequalityI}) remains such.

If $\nu_p(T)<0$, then $\nu_p(x_t)$ is initially bounded below by a constant.
If (\ref{eq:NonDegeneracy}) holds, then the minimum is achieved by a single affine
term, and (\ref{eq:InequalityI}) becomes an equality.

Given a similar analysis for $\nu_p(y_t)$, we have thus proved part $i)$ of 
theorem \ref{thm:LocalHeight}, under the restriction (\ref{eq:NonDegeneracy}) 
or the corresponding restriction for $y$ (a single restriction will suffice). 
Such a restriction avoids the pathologies described in section
\ref{section:Eigenspaces}.

\bigskip\noindent
$ii)$ \/ $\nu_p(D)<2\nu_p(T)$:
\begin{equation}\label{eq:InequalityII}
\begin{array}{llll}
\nu_p(x_t)&\geqslant&
\displaystyle \min\left\{\nu_p(x_1')+\frac{t-1}{2}\nu_p(D), 
\right.\\
&&
\displaystyle \qquad\quad \left. \nu_p(x_0')+ \frac{t-1}{2}\,\nu_p(D)
          +\nu_p\bigl(\frac{t-1}{2}\bigr)+\nu_p(T),\nu_p(x^*)\right\} 
&\quad t\, \mbox{odd}\\
\nu_p(x_t)&\geqslant&
\displaystyle  \min\left\{\nu_p(x_1')+ \frac{t-2}{2}\nu_p(D) 
          +\nu_p\bigl(\frac{t}{2}\bigr)+\nu_p(T),
\right.\\
&& 
\displaystyle \qquad\quad \left. \nu_p(x_0')+ \frac{t}{2}\,\nu_p(D),\nu_p(x^*)\right\}
&\quad t\, \mbox{even.}
\end{array}
\end{equation}
If $\nu_p(D)\geqslant 0$, then the linear terms are increasing, and we have 
two possibilities. If $x^*\not=0$, then eventually we have $\nu_p(x_t)=\nu_p(x^*)$. 
If $x^*=0$, then (\ref{eq:InequalityII}) becomes an equality provided that (here for odd $t$)
\begin{equation}\label{eq:NonDegeneracyII}
\nu_p(x_1')-\nu_p(x_0')-\nu_p(T)\not=\nu_p\bigl(\frac{t-1}{2}\bigr)
\end{equation}
and similarly for even $t$. 
The right-hand side of (\ref{eq:NonDegeneracyII}) is non-negative and grows without bounds
but at most logarithmically, due to (\ref{eq:ValuationEstimate}).
If $\nu_p(x_1')=\nu_p(x_0')$, then the left-hand side of (\ref{eq:NonDegeneracyII}) 
is negative, so this condition always holds and we have
$$
\lim_{t\to\infty} \frac{\nu_p(x_t)}{t}=\frac{\nu_p(D)}{2}.
$$
If $\nu_p(x_1')\not=\nu_p(x_0')$, then the left-hand side of (\ref{eq:NonDegeneracyII}) 
can be made negative by multiplying the initial conditions by a suitable power of $p$. 
Thus there is a rescaled sequence for which the above limit holds.
The linearity of the system ensures that the same limit holds for the original sequence.

If $\nu_p(D)<0$, then the linear terms decrease, and hence become dominant.
There is a condition analogous to (\ref{eq:NonDegeneracyII}), and we reach
an analogous result. This establishes the strict inequality in part $ii)$ 
of theorem \ref{thm:LocalHeight}.

\bigskip\noindent
$iii)$\/ $\nu_p(D)=2\nu_p(T)$:
\begin{equation}\label{eq:Summary3}
\nu_p(x_t)\geqslant \min\{
\nu_p(x_1^\prime)+(t-1)\nu_p(T),
\nu_p(x_0^\prime)+t \nu_p(T),
\nu_p(x^*)\}.\nonumber
\end{equation}
In this case we only obtain a lower bound for $\nu_p(x_t)$, and analogously for
$\nu_p(y_t)$, leading to an upper bound for $h_p$. 
One verifies that the latter agrees with the value of $h_p$ given in theorem \ref{thm:LocalHeight}
for this case.

\section{Global height}\label{section:GlobalHeight}

In this section we determine the global height (\ref{eq:LogHeight}) for the rational
points of the affine map $F$ given in (\ref{eq:Faffine}).
The dynamics of $F$ on $\R^2$ is standard \cite[section 1.2]{KatokHasselblatt}. 

Let $T,D$ and $q(x)$ be as in section \ref{section:LocalHeights}.
For a rational number $x$ we shall adopt the notation
\begin{equation}\label{eq:overunder}
x=\frac{\overline{x}}{\underline{x}}
\hskip 40pt 
\overline{x}, \underline{x}\in\mathbb{Z},\quad 
\mathrm{gcd}(\overline{x},\underline{x})=1.
\end{equation}
As before, the eigenvalues of $\mathrm{M}$ are $\alpha$ and $\beta$ 
with $|\alpha|\geqslant |\beta|$.

We consider the prime divisors of the denominators of $T$ and/or $D$, and split them
into two disjoint families:
\begin{eqnarray*}
P_1 &=& \left\lbrace p: \nu_p(\underline{D}) < 2\nu_p(\underline{T})\right\rbrace\\
P_2 &=& \left\lbrace p: \nu_p(\underline{D}) \geqslant  2 \nu_p(\underline{T}),\,
   \nu_p(\underline{D})\not=0\right\rbrace.
\end{eqnarray*}
Then we define
\begin{equation} \label{eq:ConvHeight}
h^*= \sum_{p\in P_1}\nu_p(\underline{T})\log(p)
          + \frac{1}{2}\sum_{p\in P_2}\nu_p(\underline{D})\log(p)
\end{equation}
where the sum is zero if the corresponding set of primes is empty.

\begin{theorem} \label{thm:LogH}
Let $F$ and $\mathrm{M}$ be as in (\ref{eq:Faffine}).
Then for almost all rational initial conditions $z$, the logarithmic 
height $h(z)$ defined in (\ref{eq:LogHeight}) is given by
$$
h(z)=\max(0,\log |\alpha|) + h^*
$$
where $\alpha$ is a largest eigenvalue of $\mathrm{M}$ and 
$h^*$ is as in \eqref{eq:ConvHeight}.
\end{theorem}
 
\proof
We determine the height \eqref{eq:Height} of each component $x_t$ and $y_t$ of $z_t$.
In each case, this means considering their numerator and denominator
{\em after} cancelling common factors between them, so a given prime
appears in only one of $\overline{x}_t$, $\underline{x}_t$ if it appears at all.
From \eqref{eq:Adelic}, we can write:
$$  | \underline{x}_t  |\; \prod_p\, p^{-\nu_p(\underline{x}_t)} =1,$$
where the nontrivial terms in the product correspond to the prime divisors of 
$\underline{x}_t$.

We begin with the parameter value $s=(0,0)$.
From theorem \ref{thm:LocalHeight} we have
$$
\nu_p(z_t)\,\sim\,\begin{cases} 
-t\nu_p(T)   & \mbox{if}\hskip 5pt \nu_p(D)>2\nu_p(T)\\
-t\nu_p(D)/2 & \mbox{if}\hskip 5pt \nu_p(D)\leqslant 2\nu_p(T).
\end{cases}
$$
The only primes which will contribute to the logarithmic height  of $\underline{x}_t$
are the divisors of $\underline{T}$ or $\underline{D}$. The contribution of the primes 
which divide the denominator of the initial conditions is asymptotically zero.
As a result, we have
\begin{eqnarray}
 \lim_{t\to\infty}\frac{1}{t} \log |\underline{x}_t|  
    &=&  \sum_p \lim_{t\to\infty}  \frac{\nu_p(\underline{x}_t)}{t} \log p  \nonumber \\ 
    &=&  \sum_{p\in P_1} \lim_{t\to\infty}  \frac{\nu_p(\underline{x}_t)}{t} \log p 
        +\sum_{p\in P_2} \lim_{t\to\infty}  \frac{\nu_p(\underline{x}_t)}{t} \log p \nonumber \\ 
    &=&  \sum_{p\in P_1} \lim_{t\to\infty}\frac{\nu_p({x}_t)}{t} \log p 
                           +\sum_{p\in P_2}\lim_{t\to\infty} \frac{\nu_p({x}_t)}{t} \log p  \nonumber \\ 
    &=&  -\sum_{p\in P_2}\nu_{p}({T})\log p-\frac{1}{2}\,\sum_{p\in P_1} \nu_{p}({D})\log p  \nonumber \\
    &=&  \sum_{p\in P_2}\nu_{p}(\underline{T})\log p+\frac{1}{2}\,\sum_{p\in P_1} \nu_{p}(\underline{D})\log p=h^*.\nonumber \\
\end{eqnarray}
The analogous calculation for $y_t={\overline{y}_t}/{\underline{y}_t}$
means we have established
\begin{equation} \label{eq:UnderLim}
 \lim_{t\to\infty}\frac{1}{t} \log |\underline{x}_t| 
    = \lim_{t\to\infty}\frac{1}{t} \log |\underline{y}_t| = h^*.
\end{equation}
 
Now we consider the logarithmic height of \eqref{eq:LogHeight}.  
As $\overline{x}=x\underline{x}$, we can write
\begin{eqnarray}
\lim_{t\to\infty}\frac{1}{t}\log |\overline{x}_t| &=&
\lim_{t\to\infty}\frac{1}{t} \log |x_t\; \underline{x}_t| \nonumber \\
&=& \lim_{t\to\infty}\frac{1}{t}  \log |x_t|
+ \lim_{t\to\infty}\frac{1}{t} \log |\underline{x}_t|, \label{eq:SecondPart}
\end{eqnarray}
provided the separate limits exist.
To learn about the nature of $x_t$ in the argument of the first logarithm on the right, 
we need to inject information on the archimedean dynamics of $\mathrm{M}$ on $\mathbb{Q}^2$.

We begin by assuming that $\mathrm{M}$ has diagonal Jordan form. If $|\alpha|\leqslant 1$, then 
all orbits are bounded, i.e., $|x_t|,|y_t| < C$ for some real number $C$ independent of $t$. 
We have 
\begin{equation} 
  0<|\underline{x}_t| \leqslant  H(x_t)=\max(|\underline{x}_t|,|\overline{x}_t|) \leqslant C\, |\underline{x}_t| \nonumber
\end{equation}
so that
$$
\lim_{t\to\infty}\frac{1}{t}\log H(x_{t})
= \lim_{t\to\infty}\frac{1}{t} \log |\underline{x}_t|,
$$
and similarly for $H(y_t)$ and since $\log|\alpha|\leqslant 0$, we recover \eqref{eq:ConvHeight} via \eqref{eq:UnderLim}.

If $|\alpha|>1$, then (almost) all orbits in forward time escape to infinity at the rate
$$
x_t^2+y_t^2\sim |\alpha|^{2t}(x_0^2+y_0^2).
$$
Because 
$$ 
\frac{1}{2} (x^2+y^2) \leqslant \max(|x|^2,|y|^2)  \leqslant  x^2+y^2, 
$$
it follows that 
\begin{equation} \label{eq:maxlimit}
  \lim_{t\to\infty} \frac{1}{t}\,\log \max(|x_t|,|y_t|) = \log |\alpha|.
\end{equation}
We have
\begin{eqnarray} 
\frac{1}{t}\,\log \max(|x_t|,|y_t|) &=&
\frac{1}{t}\,\max(\log|x_t|,\log |y_t|)\nonumber \\
&=& \max\left(\frac{1}{t}\, (\,\log |\overline{x}_t| -  \log |\underline{x}_t|\,),
\frac{1}{t}\, (\,\log |\overline{y}_t| -  \log |\underline{y}_t|\,)\right).\nonumber
\end{eqnarray} 
From \eqref{eq:maxlimit} and the known limits \eqref{eq:UnderLim}, we learn
$$ 
h(z_0)=\lim_{t\to\infty} \frac{1}{t}\,\log \max( |\overline{x}_t|,  |\overline{y}_t|) 
=\log(|\alpha|)+h^*
$$
as desired.

If the Jordan form of $\mathrm{M}$ is not diagonal, then $\Vert z_t\Vert$ contains an
affine term which grows sub-exponentially, and the exponential terms dominate.
If $|\alpha|=1$, then $h(z)=0$. In this case $h^*=0$ (both $P_1$ and $P_2$ are empty) and
$\log|\alpha|=0$, as desired.

It remains to consider the case $s\not=(0,0)$, corresponding to a non-zero fixed point $z^*$.
If $|\alpha|<1$, then all orbits are asymptotic to the fixed point $z^*$, so $x_t\to x^*$ and
the first term on the RHS of \eqref{eq:SecondPart} vanishes, while the case $|\alpha|\geqslant 1$ 
is dealt with by the previous analysis.
Thus, asymptotically, the logarithmic height of $\overline{x}_t$ and $\underline{x}_t$ is 
the same, similarly for $\overline{y}_t$ and $\underline{y}_t$. 
From \eqref{eq:UnderLim} we see that \eqref{eq:LogHeight} has the value $h^*$.
\endproof
 
The previous theorem shows that the logarithmic height depends only on $T$ and $D$ for 
the matrix $\mathrm{M}$ as these determine the eigenvalues.  
Thus this height is preserved by conjugacy in $\GL(2,\Q)$.
Related to $\mathrm{M}$ is 
its associated companion matrix $C$, also with rational entries:
\begin{equation} \label{mat-defs}
    C \, = \, \begin{pmatrix} T & -D \\ 1 & 0 \end{pmatrix}.
\end{equation}
It is well-known that provided $\mathrm{M}$ is not a rational multiple of the 
identity matrix, then $\mathrm{M}$ is conjugate to $C$ over $\Q$. 

\section{Piecewise affine maps}\label{section:PiecewiseAffineMaps}

We consider now two-dimensional piecewise-affine maps over the rationals,
defined as follows. Given a finite or countable set $I$ of indices,
we choose a partition of $\Q^2$ into domains $\Omega_i$, with $i\in I$.
Typically, each $\Omega_i$ will be a convex (finite or infinite) polygon.
For each $i\in I$, we choose $\mathrm{M}_i\in\mbox{\rm GL}_2(\Q)$ and $s_i\in\Q^2$,
to obtain the map $F_i:\Q^2\to \Q^2$ given by $z\mapsto \mathrm{M}_iz+s_i$.
The mapping $F$ is then defined by the rule
\begin{equation}\label{eq:FGeneral}
F:\Q^2\to\Q^2\hskip 40pt z\mapsto F_i(z),\quad z\in\Omega_i.
\end{equation}
We shall assume that the partition $\{\Omega_i\}$ is irreducible,
namely that $F$ is not differentiable on the boundaries of the
domains $\Omega_i$.

To every orbit $(z_t)$ of $F$ we associate a doubly-infinite 
sequence $\sigma=(\sigma_t)\in I^\Z$ via the rule 
\begin{equation}\label{eq:Code}
\sigma_t=i \quad \Leftrightarrow \quad z_t\in\Omega_i.
\end{equation}
The maps (\ref{eq:F}) are of the type (\ref{eq:FGeneral}), with $\Omega_i=\Delta_i\times \R$.
Their symbolic dynamics (\ref{eq:Code}) is determed by the simpler condition
$$
\sigma_t=i \quad \Leftrightarrow \quad x_t\in\Delta_i.
$$

The function $z_0\mapsto \sigma(z_0)$ is not injective, and we are interested
in the structure of the sets of points which share the same code.
The map $F$ fails to be differentiable on the set of lines and segments
$\partial \Omega$, where $\partial\Omega$ is the union of the 
boundaries of the domains $\Omega_i$.
By forming all pre-images of these lines we obtain the discontinuity set $X$ of the map:
\begin{equation}
X=\bigcup_{t\geqslant 0} F^{-t}(\partial\Omega)
\hskip 40pt\partial\Omega=\bigcup_{i\in I}\partial\Omega_i.
\end{equation}
The set $X$ is a union of segments, lines, and rays. Now consider the complement of 
the closure of $X$ in $\R^2$. This is an open set, which decomposes as the
union of connected components. By construction, all points of each connected 
component have the same code. 

The bounded connected components with a periodic code are called 
{\it islands}, denoted by $\cE$. (This terminology is normally reserved
for the area-preserving case, for which $\cE$ is also periodic.)
If $n$ is the period of the code, then $F^n$
is affine and the results of the previous section apply.  
The Jacobian $\mathrm{J}$ of $F^n$ is the same at every point of the island, since it depends 
only on the code. 
Since $\cE$ is bounded, the eigenvalues of $\mathrm{J}$ are necessarily in 
the closed unit disc in $\C$.

Let $P$ be the set of prime divisors of the denominator of the trace or 
the determinant of the matrices $\mathrm{M}_i$. 
This is the set of primes of interest to us (see also the appendix). 
Now fix $p\in P$ and embed the rational points of an island $\cE$ in the space $\Q_p^2$. 
The following result justifies the presence of plateaus in the graph of $h_p$ displayed 
in figure \ref{fig:Height}.

\begin{theorem}\label{theorem:ConstantHeight}
Almost all points of a rational island have the same heights $h$ and $h_p$ for 
all primes $p$. The latter are rational numbers.
\end{theorem}

\proof
Let $n$ be the period of the island. 
If the restriction of $F^n$ to $\cE$ has finite order, 
then all points in $\cE$ are periodic, and their height is zero. 
Let us thus assume that $F^n$ has infinite order and let $\mathrm{J}$ be 
the Jacobian of $F^n$ on $\cE$.
The result follows by applying Theorems  \ref{thm:LocalHeight} and \ref{thm:LogH}, 
respectively, to the affine map $F^n$, noting that $T$ and $D$ of
\eqref{eq:TD} now refer to the trace and determinant of $\mathrm{J}$, plus the respective 
results $h_p$ and $h^*$ of these theorems should be divided by $n$ to account for the
different time scale of the return map to the island. 
So the $p$-adic heights are rationals, in general.
\endproof

Let us now consider the behaviour of $\nu_p(z_t)$ for points in an island (figure \ref{fig:nu}).
This is the case $i)$ of theorem \ref{thm:LocalHeight}, where $\mathrm{M}=\mathrm{J}$
is the Jacobian of the return map to the island. 
The conditions of lemma \ref{lemma:SimultaneousApproximation} are 
satisfied by $\mathrm{J}$. Hence, by choosing $\zeta\in W_p^\beta$,
we can find near every point of the island initial conditions for 
orbits which perform rotations in $\Q^2$, while they simultaneously 
approach the unstable fixed point $z^*$ in $\Q_p^2$ as close as we please.

\section{Numerical experiments}\label{section:NumericalExperiments}

\begin{figure}[t]
\hfil
\epsfig{file=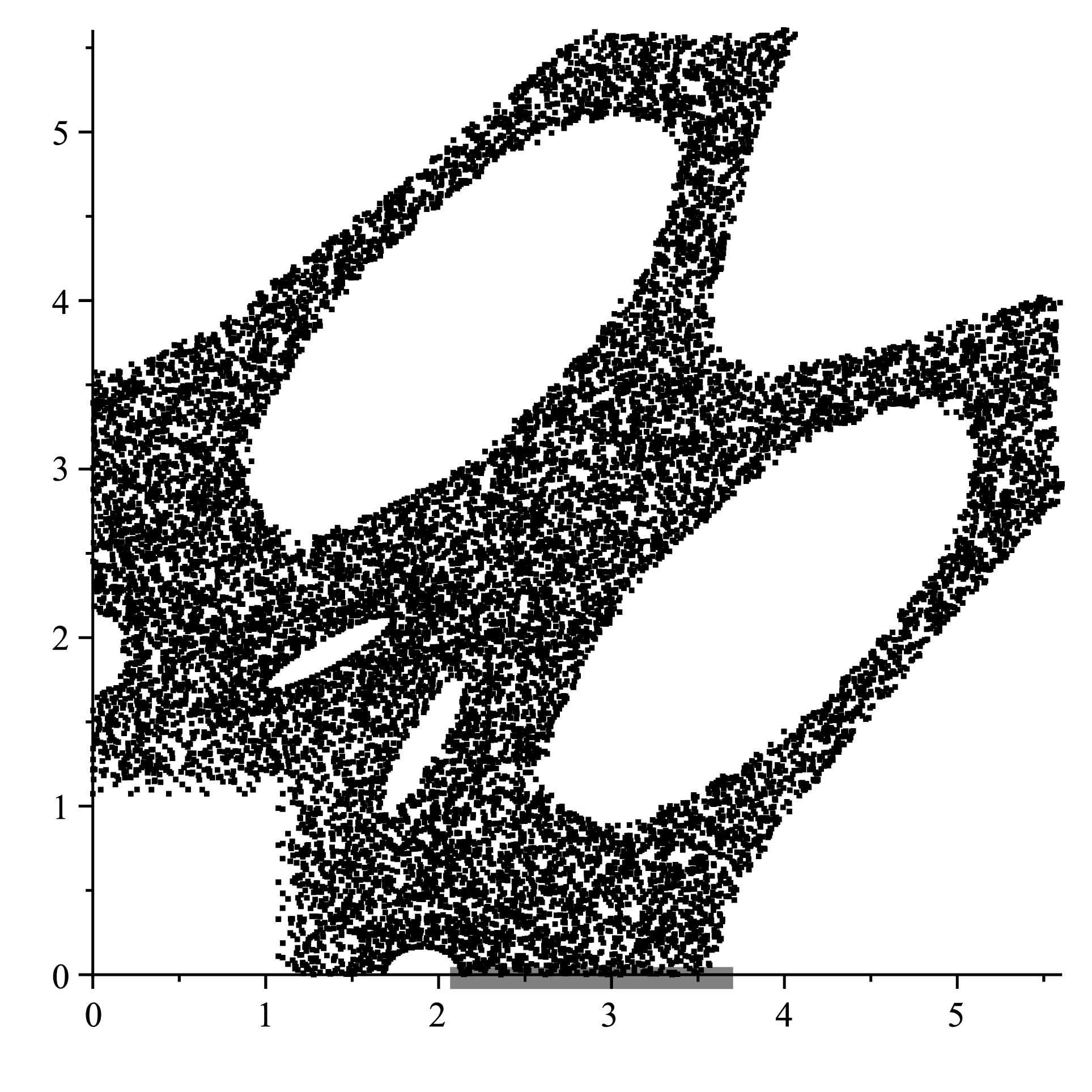,width=6cm,height=6cm}
\quad
\epsfig{file=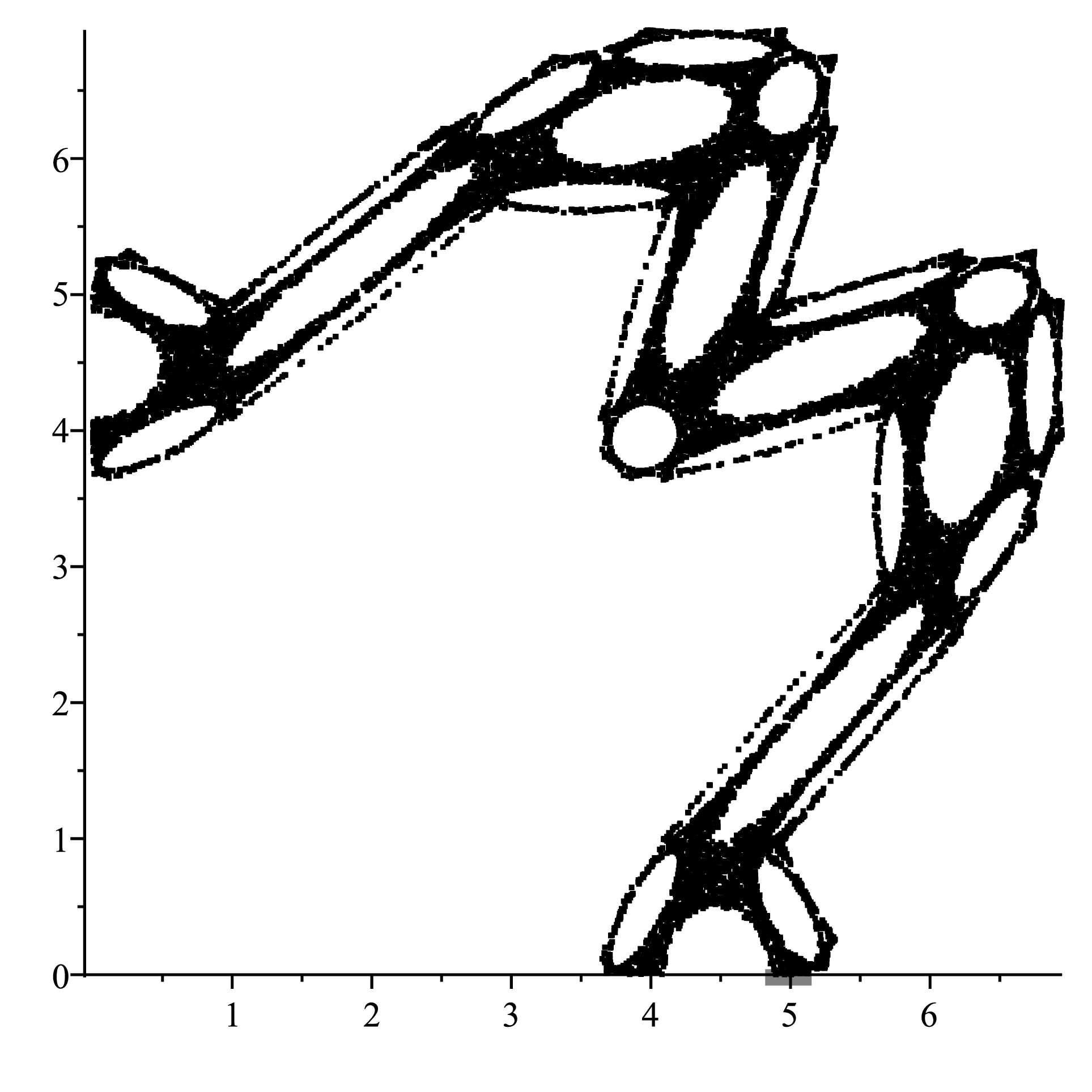,width=6cm,height=6cm}\hfil
\vspace*{-5pt}
\caption{\label{fig:ChaoticRegions}\rm\small
Chaotic regions of the map (\ref{eq:Myf}). We display the
first 50000 iterates of the point $z_0=(7/3,0)$ (left) and
$z_0=(5,0)$ (right) within the first quadrant.
}
\end{figure}

In this section we explore the convergence of heights for rational 
orbits in chaotic regions and their boundaries. Two such regions are 
displayed in figure \ref{fig:ChaoticRegions}, where in each case we 
have plotted a large number of points of a single rational orbit. 
These plots suggest that the closure of these orbits is a bounded 
subset of the plane, with positive Lebesgue measure.

At present, statements on this kind can only be established in very
special cases. For piecewise affine symplectic maps, our knowledge of the 
boundary of chaotic regions is inadequate, and proofs of global stability have
relied on the presence of piecewise-smooth bounding invariant curves, 
which is a non-generic situation \cite{Devaney,BeardonBullettRippon,LagariasRains:b}.
In the present examples there are no such curves, and we can only establish 
boundedness inside island chains. Thus any consideration on convergence of the 
height along other types of non-periodic orbits will necessarily be speculative.

\begin{figure}[h]
\hfil
\epsfig{file=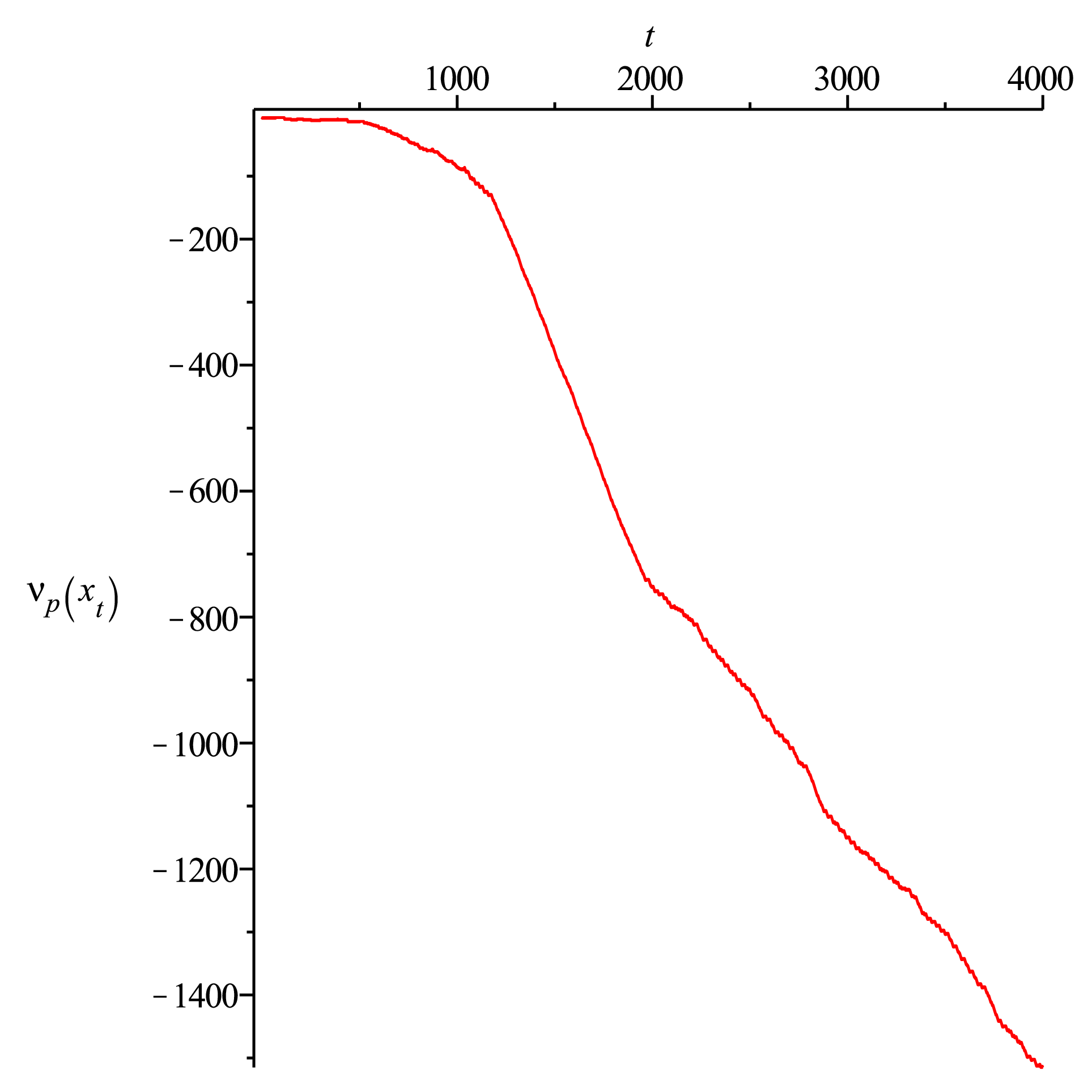,width=6cm,height=6cm}
\quad
\epsfig{file=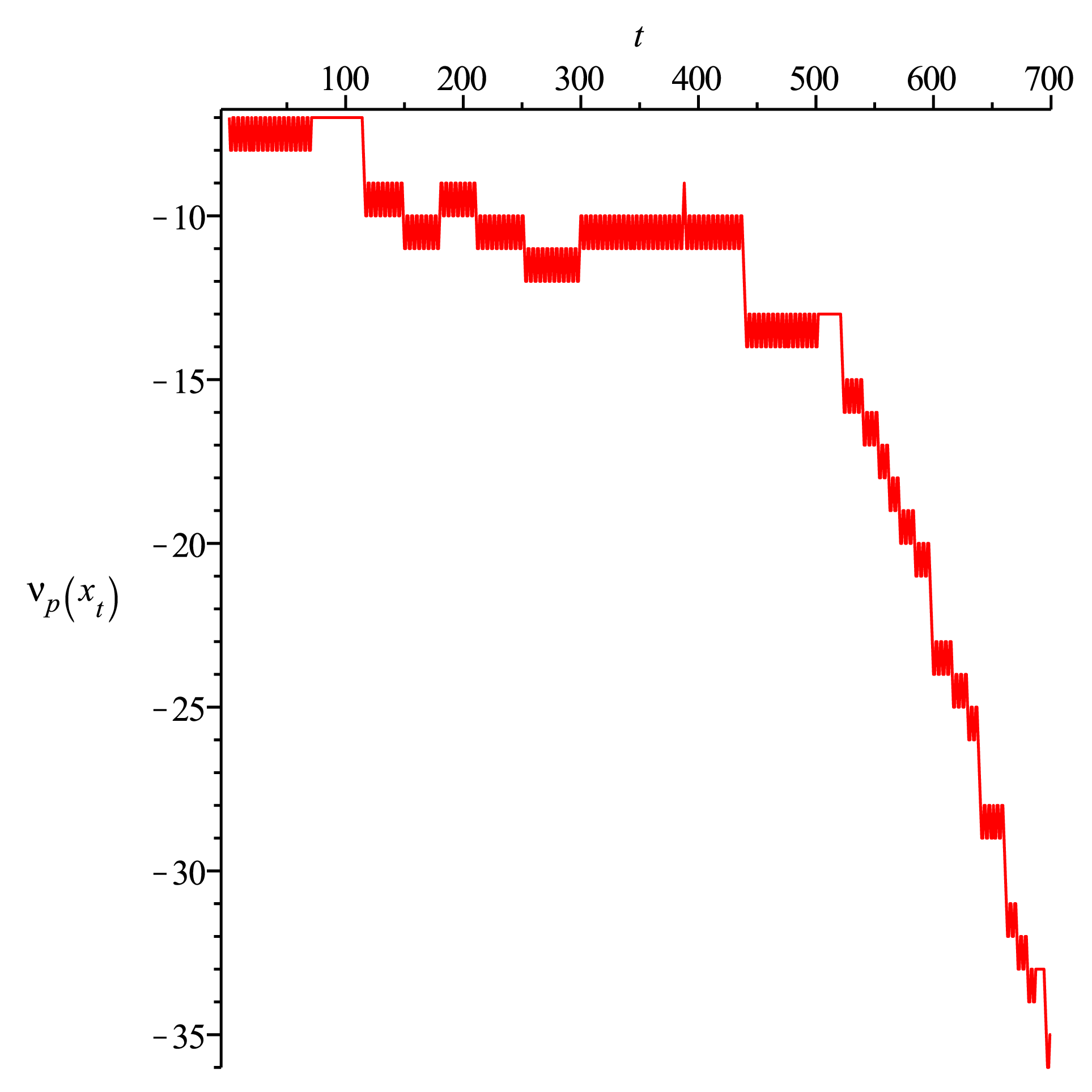,width=6cm,height=6cm}\hfil
\vspace*{-5pt}
\caption{\label{fig:ValuationVsTime}\rm\small
Time-dependence of $\nu_2(x_t)$ for one rational orbit 
of the map (\ref{eq:Myf}).
Left: typical behaviour, showing transitions between four 
different regimes. Right: detail of the first plateau and 
the beginning of the drop.
}
\end{figure}

We begin to examine the behaviour of $\nu_p(x_t)$ along an individual orbit 
of the area-preserving map (\ref{eq:F}), with $f$ given by (\ref{eq:Myf}).
There is only one prime in $P$, namely $p=2$ (the set $P$ was defined in section
\ref{section:PiecewiseAffineMaps}). 
We choose the initial condition $z_0$ near the boundary 
of the square stable region containing the origin in figure \ref{fig:ChaoticRegions}, left.
The time-dependence of $\nu_2$, shown in figure \ref{fig:ValuationVsTime},
features a concatenation of distinct regimes, in which the rate of change 
of $\nu_2$ remains approximately constant.

\begin{figure}[t]
\hfil
\epsfig{file=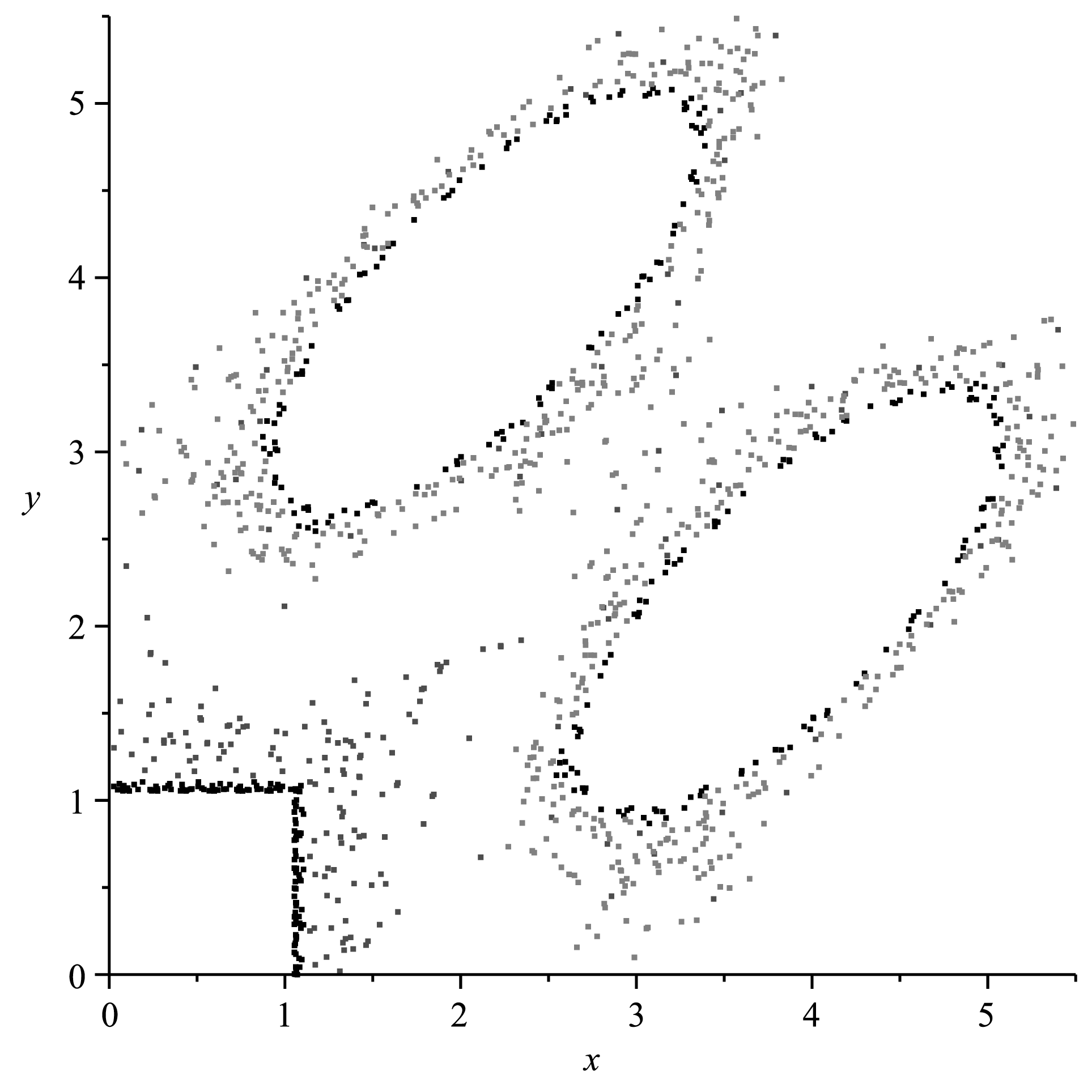,width=7cm,height=7cm}\hfil
\vspace*{-5pt}
\caption{\label{fig:ChaoticOrbit}\rm\small
Phase plot of the orbit of figure \ref{fig:ValuationVsTime}. 
The points corresponding to the four different sections
of the left diagram are plotted in different shades of grey.
}
\end{figure}

Each regime has a dynamical signature.
In figure \ref{fig:ChaoticOrbit} we plot the orbit that
generated the data of figure \ref{fig:ValuationVsTime}.
The initial plateau corresponds to the neighbourhood
of the square island mentioned above. After a transitional 
phase, the orbit migrates to a neighbourhood of the large 
island chain visible in the middle of the chaotic sea, 
where the local value of the height (the slope of the curve) 
remains approximately constant. 
Then the orbit leaves this region, and the height
decreases.

\begin{figure}[h]
\hfil
\epsfig{file=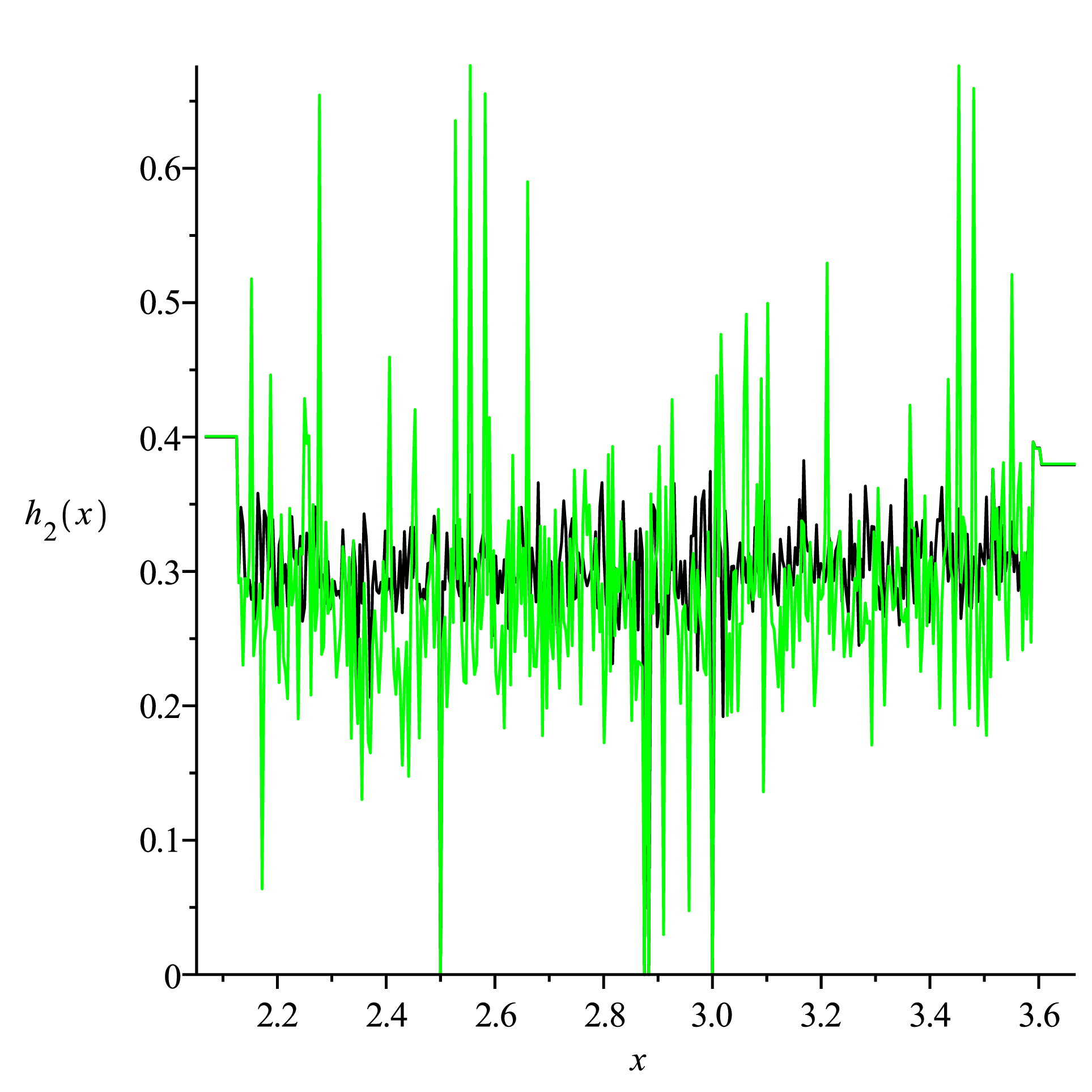,width=6cm,height=6cm}
\quad
\epsfig{file=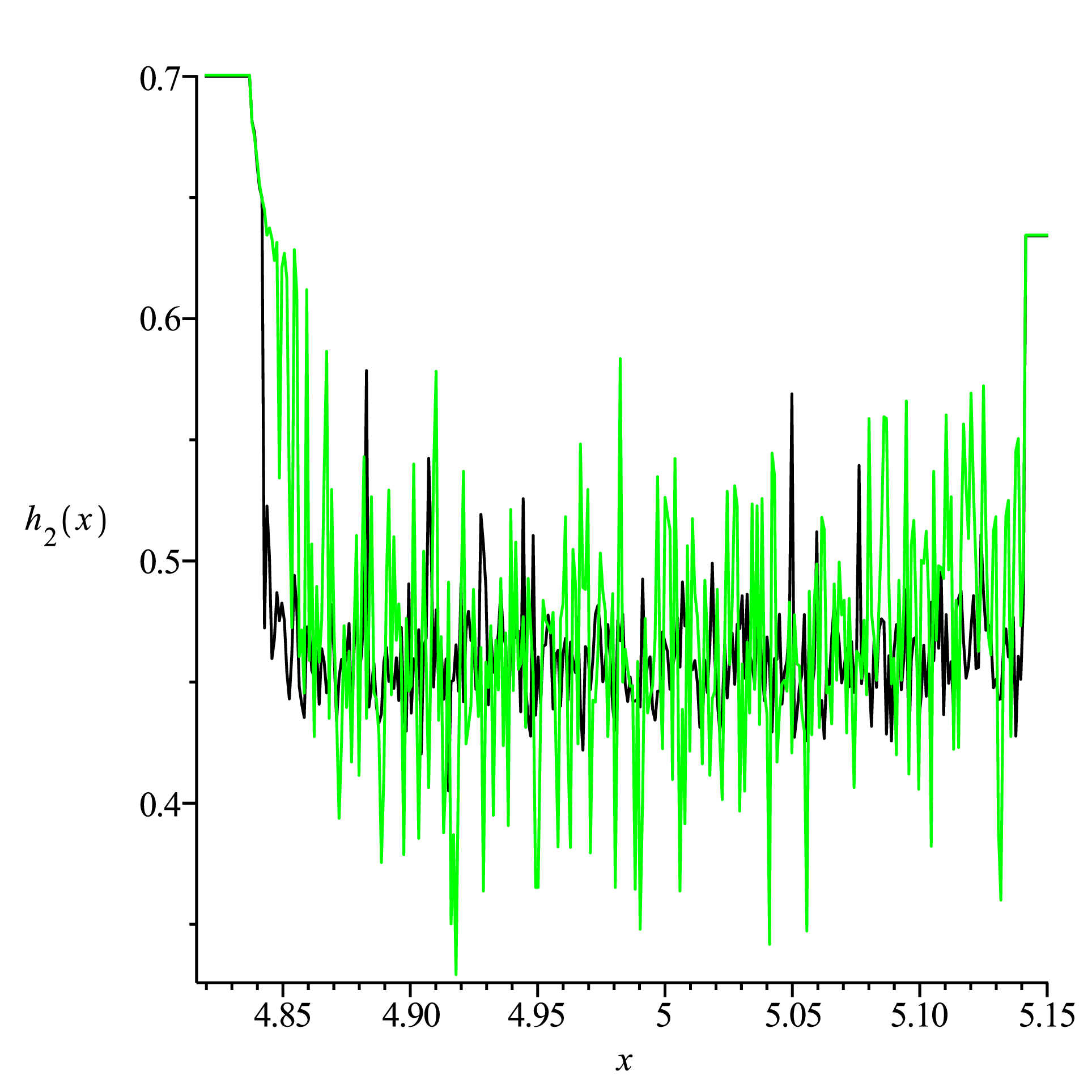,width=6cm,height=6cm}\hfil
\vspace*{-5pt}
\caption{\label{fig:HeightChaoticRegion}\rm\small
Value of $-\nu_2(x_T)/T$ for approximately 300 orbits with initial 
conditions $z_0^{(i)}, i=1,2,\ldots$ evenly spaced along a segment 
in phase space. The end-points of the segment lie inside islands,
where the height is constant. In both figures the black and green 
curves correxpond to to $T=4000$ and $T=64000$, respectively, indicating slowly
decreasing fluctuations.
}
\end{figure}

To shed light on the global picture, we have computed the approximate 
value of the height for some 300 distinct orbits. The initial 
conditions are points equally spaced on a segment connecting 
two islands, but otherwise lying in the chaotic sea. These segments 
are placed along the $x$-axis, and are visualized as grey strips 
in figure \ref{fig:ChaoticRegions}.
A numerical approximation for the height, given by 
\begin{equation}\label{eq:Approximation}
h_p(z_0,T)= \frac{\nu_p(z_0)-\nu_p(z_T)}{T} \approx h_p(z_0)
\end{equation} 
is computed for each orbit at several values of $T$: $T=4000,8000,16000,32000,64000$. 
(The value $T=64000$ yields rational
numbers with over 5000 decimal digits at numerator and denominator.) 
The data for $T=4000$
and $T=64000$ are displayed in figure \ref{fig:HeightChaoticRegion}.

\begin{figure}[h]
\hfil
\epsfig{file=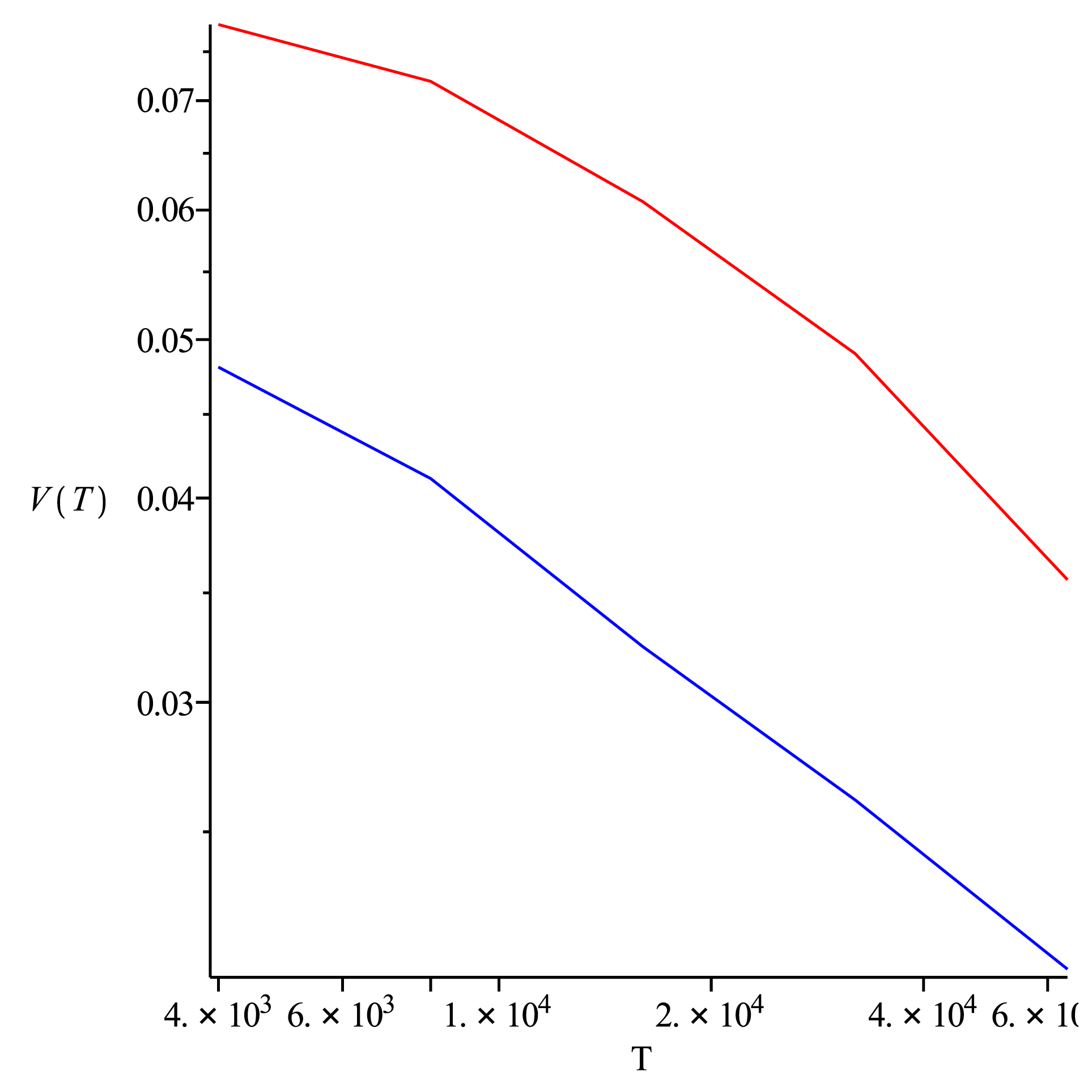,width=8cm,height=8cm}\hfil
\vspace*{-5pt}
\caption{\label{fig:VariationHeightChaoticRegion}\rm\small
Plot of $V_N(T)$ defined in (\ref{eq:VariationHeightChaoticRegion})
versus the number $T$ of iterations for the data of figure
\ref{fig:HeightChaoticRegion} (the red and blue curves correspond
to the left and right plots in the figure, respectively).
}
\end{figure}

The fluctuations appear to decrease, albeit slowly, with $T$.
To quantify this phenomenon we have computed the normalised
total variation $V$ of the height
\begin{equation}\label{eq:VariationHeightChaoticRegion}
V_N(T)=\frac{1}{N-1}\sum_{i=1}^{N-1}\bigl|h_2(z_0^{(i+1)},T)-h_2(z_0^{(i)},T)\bigr|
\end{equation}
where $N$ is the number of orbits, and $z_0^{(i)}$ is the
initial condition of the $i$th orbit.
The behaviour of $V_N(T)$ for both cases is shown in figure
\ref{fig:VariationHeightChaoticRegion} in doubly logarithmic scale. 
The data suggest a regular decrease of the total variation of the numerical 
height with the time $T$, and are consistent with a slow convergence to 
a value which is constant almost everywhere in a chaotic region.
Clearly there will be exceptional orbits where the height assumes a different
value, such as unstable periodic orbits.

\begin{figure}[h]
\epsfig{file=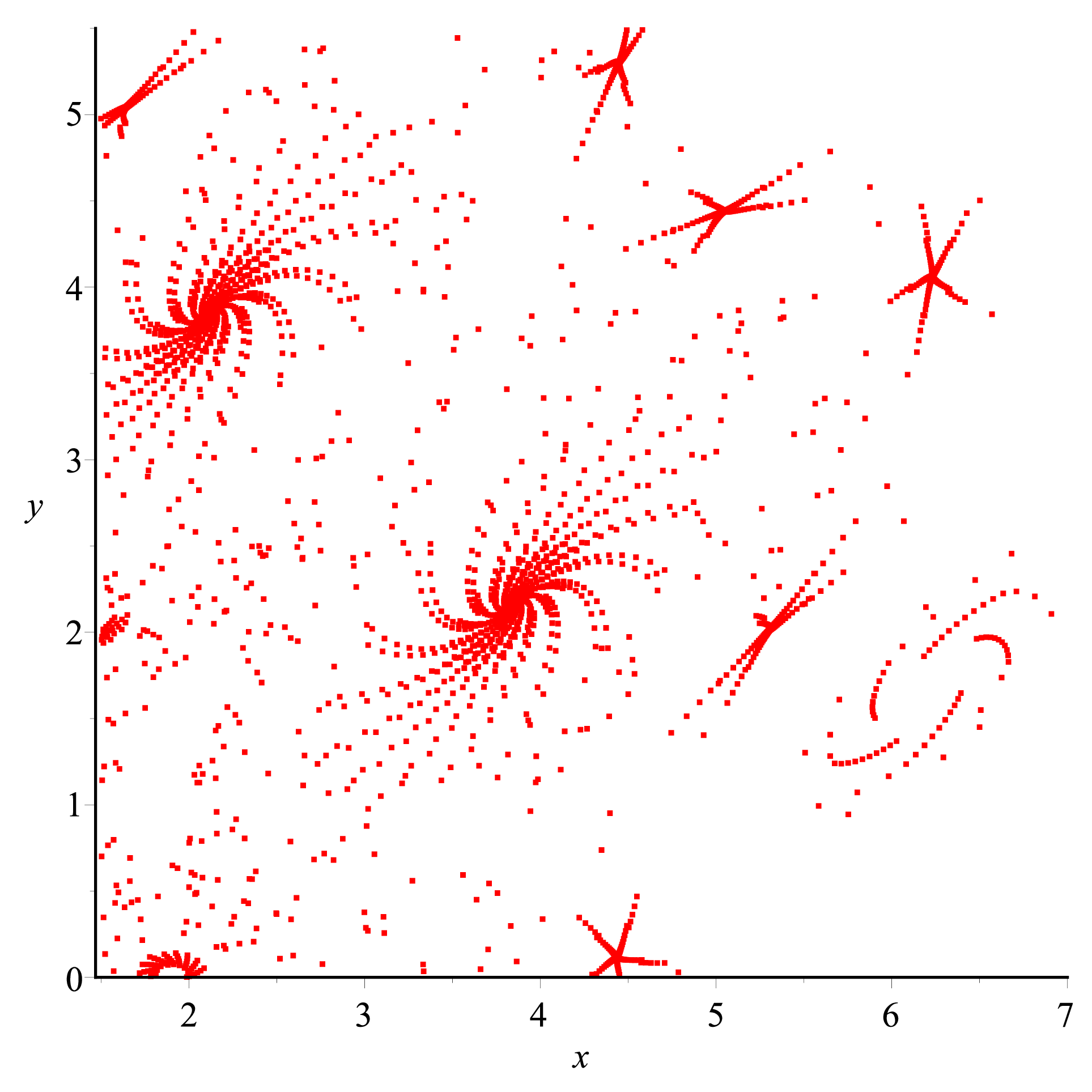,width=7cm,height=7cm}
\quad
\epsfig{file=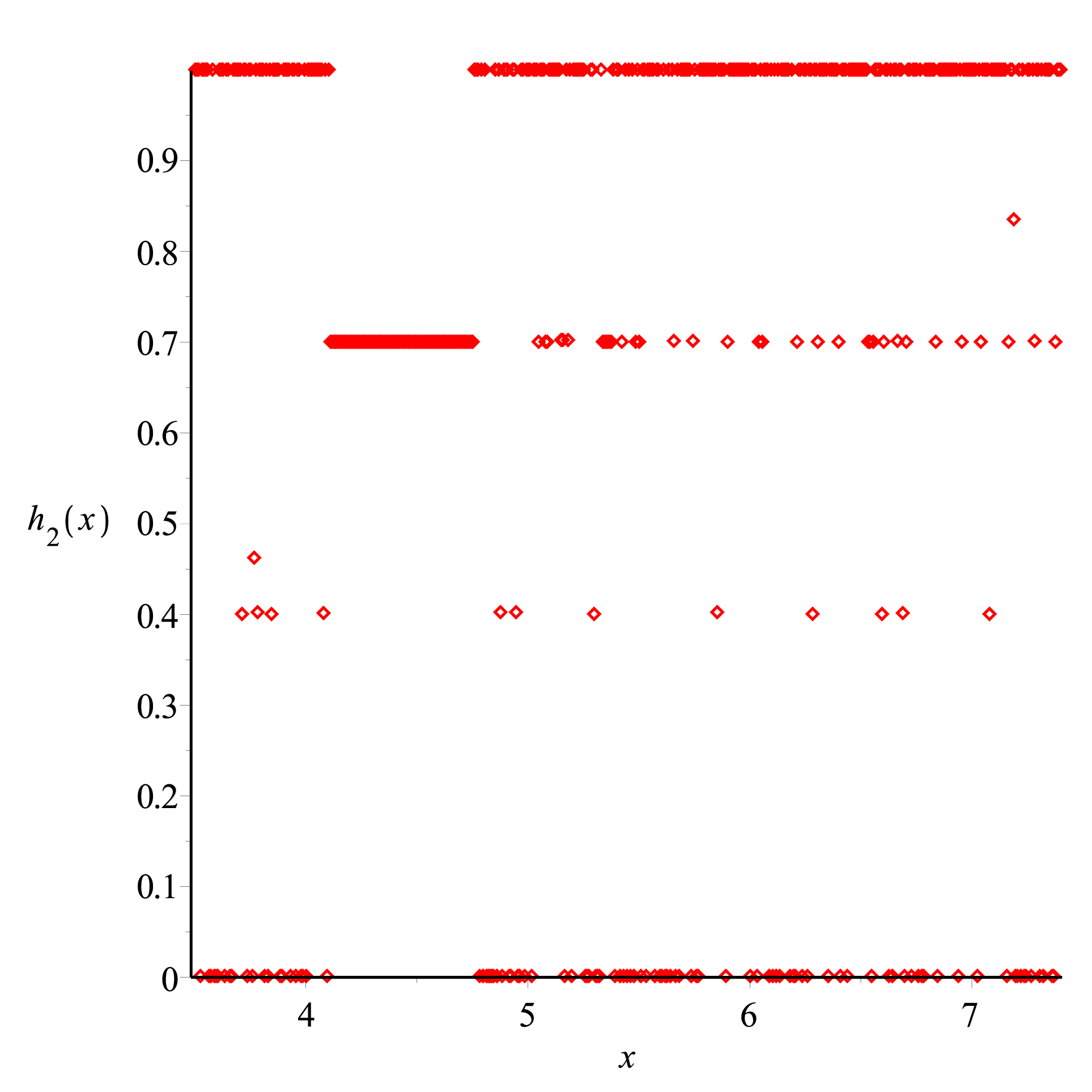,width=7cm,height=7cm}
\vspace*{-5pt}
\caption{\label{fig:PhasePortraitDiss}\rm\small
The dissipative map $F$ given in (\ref{eq:F}), with $f$ as in
(\ref{eq:Myf}) (the same as in figure \ref{fig:PhasePortrait})
and $d=497/499$.
Left: Phase portrait, with orbits spiralling towards the centres
of the islands.
Right: the $2$-adic height $h_2(x)$ for $z_0=(x,0)$ (to be compared
with figure \ref{fig:Height}, right). The limited set of values it
assumes (four, in total) reflects the existence of a limited
number of attractors. The absence of fluctuations indicates
that these attractors have a simple structure.
}
\end{figure}

The scenario for dissipative maps is simpler; the orbits, after a transient,
relax to a small number of point attractors (figure \ref{fig:PhasePortraitDiss}, left).
In figure \ref{fig:PhasePortraitDiss}, right, we plot the approximate height 
$h_2$ for initial conditions of the type $z_0=(x,0)$ with $x$ in an interval 
which crosses an island. Outside the islands the height jumps wildly between few 
values, presumably due to the very complicated boundaries of the basins of the
various attractors.

We synthesise our findings with two conjectures. 

\medskip
\noindent{\bf Conjecture 1}. {\sl Let $f$ be a piecewise affine map of 
$\mathbb{Q}^2$ and let $O$ be a bounded orbit of $f$. 
Then, for each prime $p$, the functions $h$ and $h_p$ are almost 
everywhere constant on $\overline O\cap \mathbb{Q}^2$, where 
$\overline O$ is the closure of $O$ in $\R^2$. 
}

\medskip
Here the term `almost everywhere' refers to full density in expression (\ref{eq:Density}).

\medskip
\noindent{\bf Conjecture 2}. {\sl Let $f$ and $O$ be as above, and let
$O$ have zero Lyapunov exponent. Then, for any $p\in P$, the height 
$h_p, \, p\in P$ has a (non-strict) local maximum at $O$.
}

\begin{figure}[t]
\hfil
\epsfig{file=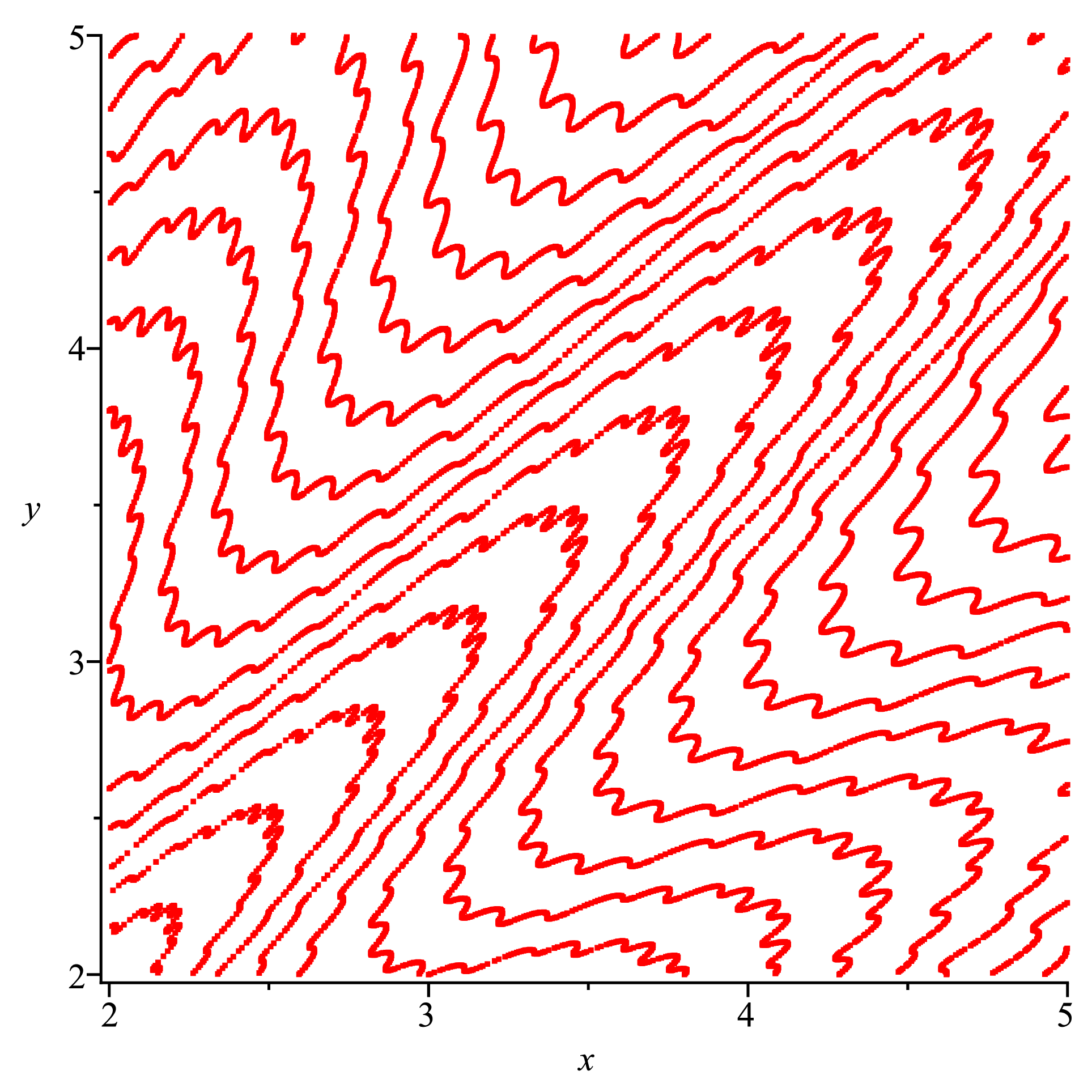,width=7cm,height=7cm}\hfil
\vspace*{-5pt}
\caption{\label{fig:Lagarias}\rm\small
Foliation of the plane into non-smooth invariant curves 
for the map (\ref{eq:Lagarias}) for $a_1=2/3, a_2=3/2$.
}
\end{figure}

\medskip
In the present context, we have identified regular orbits with linear bounded 
orbits within islands, which either foliate the island into invariant ellipses or
spiral towards the fixed point in the centre.
No analysis of planar maps would be complete without some reference to more 
general types of regular orbits, namely quasi-periodic orbits on invariant curves 
(not necessarily smooth) which are topologically conjugate to 
irrational rotations.
It has long been known that non-smooth symplectic maps may support isolated 
invariant curves \cite{HenonWisdom}, and even foliations of non-smooth curves,
see figure \ref{fig:Lagarias}.
Unfortunately the existence of such curves ---isolated or not--- for non-smooth maps 
cannot be established in general, and this limitation applies to piecewise affine maps 
with rational parameters considered here.

There are however important results for specific models.
These include a specific two-parameter family of piecewise-linear 
mappings of the type (\ref{eq:F}), where a foliation of the plane into invariant
curves has been proved (or can reasonably be conjectured) to exist
\cite{BeardonBullettRippon,LagariasRains,LagariasRains:b,LagariasRains:c}.
These are maps or type (\ref{eq:F}), with the piecewise linear functions
\begin{equation}\label{eq:Lagarias}
f(x)=\begin{cases}
 a_1x & x < 0\\
 a_2 x & x \geqslant 0.
\end{cases}
\end{equation}
Due to local linearity, these maps transform the lines through the origin into 
themselves while preserving their order, thereby inducing a circle map with a 
well-defined rotation number. 

The existence of piecewise-smooth invariant curves has been established 
for some parameter values given by algebraic numbers of degree 2 
\cite[theorem 2.2]{LagariasRains:b}.
The situation for rational parameters less clear.
If the rotation number is irrational with bounded partial quotients, then 
an early result by M.~Herman \cite[theorem VIII.5.1]{Herman} implies that
(\ref{eq:Lagarias}) is topologically conjugate to a planar rotation. 
To the authors' knowledge, the required diophantine condition have not been
established in the case of rational parameter $a_{1}$ and $a_2$ in 
(\ref{eq:Lagarias}).

Numerical experiment suggest that for rational parameters $a_1\not=a_2$, 
if the orbits of the map $f$ are bounded, then the plane foliates into invariant 
curves which typically are non-smooth. 
Under such circumstance, we found that all height functions are constant over 
the entire plane. 
This suggests that conjectures 1 and 2 hold for orbits on invariant curves as well.

\section*{Appendix}\label{section:Appendix}

We define a module $\mathbb{L}$ with the property that $\mathbb{L}^2$ serves as a minimal phase 
space for piecewise-affine maps $F$ of the form (\ref{eq:FGeneral}) with 
$F_i(z)=\mathrm{M}_iz+s_i$.
Let $\mathrm{M}_i=(m^{}_{j,k})$ and let $P$ be the (possibly empty, or infinite) set of primes which 
divide the denominator of $m^{}_{j,k}$ for some $j,k$. 
If $P$ is empty, then we let $\K=\Z$; otherwise we let
\begin{equation}\label{eq:Ring}
\K=\prod_{p\in P}\Z\left[\frac{1}{p}\right]
\end{equation}
where the product denotes the algebraic (Minkowski) product of sets.
The set $\K$ is the sub-ring of $\Q$ consisting of all the rationals 
whose denominator is divisible only by primes in $P$. 
The module $\mathbb{L}$ of the map $F$ is defined as 
\begin{equation}\label{eq:Module}
\mathbb{L}=\K+\sum_{{i\in I}\atop{j=1,2}}\{s^{(i)}_j\}
\end{equation}
where $s_i=(s^{(i)}_1,s^{(i)}_2)$ and the sum denotes algebraic sum of sets.
The set $\mathbb{L}$ is a $\K$-module (a group under addition, with 
a multiplication by elements of $\K$). 

If $I$ is finite, then there is an integer $N$ such that
$$
\mathbb{L}=\frac{1}{N}\,\K.
$$
To compute $N$, we let $d_i$ be the least common multiple of the denominators
of $s^{(i)}_1$ and $s^{(i)}_2$  and let
\begin{equation}\label{eq:dprime}
d_i'=d_i\prod_{p\in P}p^{-\nu_p(d_i)}\qquad i\in I.
\end{equation}
(This product is finite.) Thus $d_i'$ is the largest divisor of $d_i$ 
which is co-prime to all primes in $P$. Then $N$ is the least common
multiple of the $d_i'$s, for $i\in I$.

If $I$ is infinite, then the integer $N$ defined above need not exist.

By construction, we have that $F_i(\mathbb{L}^2)\subset \mathbb{L}^2$ for all $i\in I$.
Hence $F(\mathbb{L}^2)\subset \mathbb{L}^2$ and $\mathbb{L}^2$ is a natural minimal phase space
for $F$. 

The set $\mathbb{L}$ may be embedded in $\Q_p$ for any prime $p$ 
(the field $\Q_p$ is the completion of $\Q$ with respect to the absolute value $|\cdot|_p$).
If $p\in P$, then $\mathbb{L}$ is an unbounded dense subset, and so even if the $\Q^2$ 
motion is bounded, the $p$-adic dynamics may be unbounded.
If $p\not\in P$ and the set $I$ of indices is finite, then $\mathbb{L}$ 
is bounded in $\Q_p$, and if $p$ does not divide any of the $d_i'$ 
(see (\ref{eq:dprime})), then $\mathbb{L}$ lies within the unit disc in $\Q_p$. 
If $I$ is infinite, then $\mathbb{L}$ may still be unbounded even if $p\not\in P$, 
that is, the $p$-adic height may grow entirely due to the additive action 
of $F$ (the translations $s_i$).


\end{document}